\documentclass[11pt,a4paper]{article}
\usepackage[T1]{fontenc}
\usepackage[utf8]{inputenc}
\usepackage{amsmath,amssymb,amsthm,mathtools}
\usepackage{newtxtext,newtxmath}
\usepackage{microtype}
\usepackage[a4paper,top=25mm,bottom=25mm,left=25mm,right=25mm,headheight=14pt]{geometry}
\usepackage{authblk}
\usepackage{enumitem}
\usepackage{array,booktabs}
\usepackage{placeins}
\usepackage{needspace}
\usepackage{fancyhdr}
\usepackage[hidelinks]{hyperref}
\usepackage[nameinlink,noabbrev]{cleveref}
\numberwithin{equation}{section}
\newtheorem{theorem}{Theorem}[section]
\newtheorem{proposition}[theorem]{Proposition}
\newtheorem{lemma}[theorem]{Lemma}
\newtheorem{corollary}[theorem]{Corollary}
\theoremstyle{definition}
\newtheorem{definition}[theorem]{Definition}
\newtheorem{example}[theorem]{Example}
\theoremstyle{remark}
\newtheorem{remark}[theorem]{Remark}
\newcommand{\RR}{\mathbb R}
\newcommand{\ZZ}{\mathbb Z}
\newcommand{\TT}{\mathbb T}
\newcommand{\PP}{\mathbb P}
\newcommand{\EE}{\mathbb E}
\newcommand{\cG}{\mathcal G}
\newcommand{\cS}{\mathcal S}
\newcommand{\cY}{\mathcal Y}
\newcommand{\cP}{\mathcal P}
\newcommand{\dd}{\,\mathrm d}
\newcommand{\one}{\mathbf 1}
\newcommand{\SVU}{\ensuremath{\mathrm{SVU}}}
\newcommand{\W}{\mathsf W_1}
\newcommand{\Wsyn}{\mathsf W_{1,\mathrm{syn}}}
\DeclareMathOperator{\Prob}{Prob}
\DeclareMathOperator{\supp}{supp}
\DeclareMathOperator{\Gr}{Gr}
\DeclareMathOperator{\SL}{SL}
\DeclareMathOperator{\GL}{GL}
\DeclareMathOperator{\dist}{dist}
\DeclareMathOperator*{\esssup}{ess\,sup}
\setlist[enumerate]{label=\textup{(\roman*)},leftmargin=*,itemsep=3pt,topsep=4pt}
\fancypagestyle{plain}{%
  \fancyhf{}\fancyfoot[C]{\thepage}%
}

\title{\vspace{-1.7em}\bfseries Exponential Deviations and H\"older Continuity\\[0.2em]
for Moving-Law Mixed Cocycles}
\author{Zhongmin Liu}
\author{Dengdi Chen}
\author{Yan Zheng\thanks{Corresponding author: Yan Zheng.}}
\affil{College of Science, National University of Defense Technology, Changsha 410073, PR China\\
\texttt{1254642090@qq.com} \quad \texttt{dengdichen@yeah.net} \quad \texttt{yanzhengyl@163.com}}
\date{}
\begin{document}
\maketitle\vspace{-2.0em}
\begin{abstract}
We study mixed random–quasiperiodic cocycles whose one-step laws evolve along an ergodic environment. Under mixed-base ergodicity, bounded logarithmic growth, and stationary-value uniqueness, which we characterize by relative quasi-irreducibility, we prove exponential quenched deviations for finite-time directional growth.
A renewal example shows that, in general ergodic environments, annealed exponential deviations can fail despite quenched exponential concentration. This obstruction motivates us to consider a fixed finite-state irreducible aperiodic Markov environment, where we establish annealed exponential deviation estimates uniform over a relative neighborhood of a reference sampling kernel in the bounded mixed-base-ergodic class at each fixed tolerance.
These estimates yield local Hölder continuity of the maximal Lyapunov exponent near a stationary-value-unique kernel with a simple top exponent. Under common matrix and inverse bounds, the continuity applies to synchronized Wasserstein perturbations with fixed translation marginals and continuous matrix fields, or to ordinary Wasserstein perturbations with uniform Hölder regularity in the phase.
As an application, suitable frequency-support and noise-nondegeneracy conditions yield Hölder energy dependence for moving-law Schrödinger cocycles with continuous potentials.
\end{abstract}
\vspace{0.6em}
\noindent\textbf{2010 Mathematics Subject Classification.}
Primary 37H15; Secondary 37A25, 60B20, 60G42.\par
\vspace{0.3em}
\noindent\textbf{Keywords.}
random cocycles, ergodic environments, stationary projective families,
exponential deviation estimates, Lyapunov exponents, H\"older continuity.\par
\vspace{0.3em}
\noindent\textit{Funding.} Supported by the National Natural Science Foundation of China under Grant 12571189.
\vspace{1.4em}

\section{Introduction}\label{secintro}

Lyapunov exponents describe the asymptotic growth of matrix products. Two quantitative
questions arise naturally. How closely does the growth of a finite matrix product reflect
its asymptotic rate, and how stable is that rate under changes in the one-step sampling
distributions? These questions are central to random cocycles and to the transfer-matrix
analysis of Schr\"odinger operators.
A mixed random--quasiperiodic cocycle combines random translations of a torus phase with
phase-dependent linear maps on vector fibers. We study finite-time deviations and Lyapunov
regularity when an ergodic environment determines the distribution of each sampled
transformation.

The classical theories of Furstenberg, Furstenberg--Kifer, Kingman, and Oseledets give
the asymptotic foundations for random matrix products \cite{Furstenberg,FK,Kingman,Oseledets}.
Cai, Duarte, and Klein introduced mixed random--quasiperiodic cocycles with a fixed
one-step law and established Furstenberg-type positivity and expected-growth results
for that setting \cite{CDKmodel,CDKfurstenberg}. Their work on mixing Markov systems
also provides effective deviation estimates for suitable observables \cite{CDKstatistics}.
A different line of research allows the one-step distribution to depend on the environment.
Kifer studied random bundle maps whose distributions evolve with a stationary
random environment \cite{Kifer}, while Goldsheid established exponential
growth results for nonstationary Markov-dependent matrix products
\cite{Goldsheid2022}.

We formulate these questions for moving-law mixed random--quasiperiodic
cocycles driven by an invertible ergodic environment
$(\Omega,\mathcal F,P,\Theta)$. The one-step transformations belong to the
group $\cG=\TT^d\times C^0(\TT^d,\SL_r(\RR))$ of quasiperiodic
cocycle elements. An element $g=(\alpha,B)\in\cG$ acts on a phase and vector by
$(\theta,v)\mapsto(\theta+\alpha,B(\theta)v)$. At environmental state $\omega$,
the one-step element is sampled from a probability law $\nu_\omega$ on $\cG$.
Together, these one-step laws define a probability kernel
$\nu^\bullet:\omega\mapsto\nu_\omega$ from $\Omega$ to $\cG$.
Conditional on an initial environment $\omega$, let
$\xi=(\xi_j)_{j\ge0}\in\cG^{\mathbb N}$
be the sequence of sampled elements. Its coordinates are independent, with
$\xi_j$ distributed according to $\nu_{\Theta^j\omega}$. Their conditional law on
sample sequences is $\PP_\omega=\bigotimes_{j\ge0}\nu_{\Theta^j\omega}$. We
write $A_n(\xi,\theta)$ for the matrix product formed by the first $n$ elements,
with each matrix evaluated at the phase reached by the preceding translations.

For an initial phase and projective direction
$x=(\theta,\hat p)\in Q:=\TT^d\times\PP(\RR^r)$, let $p\ne0$ be a
representative of the line $\hat p$. The accumulated logarithmic growth
in this direction over $n$ steps is the random variable
\[
 Z_n^\omega(x)(\xi)
 :=\log\frac{\|A_n(\xi,\theta)p\|}{\|p\|},
 \qquad \xi\sim\PP_\omega.
\]
Under mixed-base ergodicity and logarithmic integrability, $L_1$ denotes the top
Lyapunov exponent, the almost-sure limit of $n^{-1}\log\|A_n(\xi,\theta)\|$.

Unlike the fixed-law setting, stationary projective measures in the
moving-law model are represented by environment-indexed families.
Such a family $\eta^\bullet=(\eta_\omega)$ is transported between
successive environments according to
\[
 \eta_{\Theta\omega}=\int_{\cG}g_*\eta_\omega\,\dd\nu_\omega(g).
\]
The corresponding Furstenberg formula in the moving-law setting was established in
\cite{LCZpositivity}. A key structural condition in our quantitative analysis is
\emph{stationary-value uniqueness} (\SVU), which means that all stationary
projective families have the same Furstenberg growth value $L_1$.
Under mixed-base ergodicity and logarithmic integrability, we prove that
\SVU\ is equivalent to a quasi-irreducibility condition relative to the
random environment
(see \Cref{thm:QI-SVU-equivalence}).

We first establish uniform convergence of normalized conditional expected growth
over initial phases and directions. With the environment fixed, let $\EE_\omega$
denote expectation under the conditional path law $\PP_\omega$. Set
$H_n(\omega,x):=\EE_\omega Z_n^\omega(x)$ and
\[
 \Delta_n(\omega)
 :=\sup_{x\in Q}\left|\frac{H_n(\omega,x)}{n}-L_1\right|.
\]
A compactness argument for Young measures identifies the limiting
upper and lower conditional growth envelopes with the extremal
stationary projective values. Thus we prove that \SVU\ is equivalent to $\Delta_n\to0$ in $L^1(P)$,
and also implies almost-sure convergence, as established in
\cref{thmexpected}. Compared with the fixed-law result of Cai, Duarte, and Klein
\cite{CDKfurstenberg}, our theorem dispenses with the gap between the two
largest Lyapunov exponents. Compared with Kifer's random-environment
result \cite{Kifer}, it requires neither this gap nor exponential mixing
of the underlying Markov transitions.

This convergence of conditional expected growth leads us to study
\emph{quenched} exponential deviations. We prove in \cref{thmquenched} that,
under mixed-base ergodicity, \SVU, and uniformly bounded one-step logarithmic
growth, for every $\varepsilon>0$ there are deterministic
$C_\varepsilon,c_\varepsilon>0$ and an almost surely finite measurable
$N_\varepsilon(\omega)$ such that
\begin{equation}\label{eq:intro-quenched}
 q_n(\omega,\varepsilon)
 :=\sup_{x\in Q}\PP_\omega
 \left\{\left|\frac{Z_n^\omega(x)}{n}-L_1\right|>\varepsilon\right\}
 \le C_\varepsilon e^{-c_\varepsilon n}
 \qquad(n\ge N_\varepsilon(\omega)).
\end{equation}
Unlike exponential deviation results for mixing Markov systems
\cite{CDKstatistics}, our result requires no quantitative mixing of the environment.

The \emph{annealed} deviation quantity is the environmental average
$\int_\Omega q_n(\omega,\varepsilon)\,\dd P(\omega)$.
The quenched bound \eqref{eq:intro-quenched} does not in general imply
exponential decay of this average. We show in
\cref{proprenewal-obstruction} that the annealed deviation quantity for a renewal example
is bounded below by $cn^{-2}$ at the fixed tolerance
$\varepsilon=L_1/2>0$, although it satisfies the hypotheses of
\cref{thmquenched}. Thus ergodicity of the environment alone does not
ensure exponential annealed concentration.

To obtain exponential annealed bounds, we restrict the environment to
the stationary two-sided shift of a fixed irreducible aperiodic Markov
chain with finite state space $S$. Here $\boldsymbol\nu=(\nu_s)_{s\in S}$ denotes
the sampling kernel from $S$ to $\cG$, with component one-step laws
$\nu_s$. Fix a reference kernel $\boldsymbol\nu_0$ whose mixed base is
ergodic and which satisfies \SVU. Keeping the Markov transition matrix
fixed, we vary the sampling kernel within the bounded mixed-base-ergodic
class, where the matrix fields and their inverses obey a common bound. We prove in
\cref{thmuniform-ldt} that, for each $\varepsilon>0$, there are a
neighborhood $\mathcal U_\varepsilon$ of $\boldsymbol\nu_0$ relative to
this bounded mixed-base-ergodic class, deterministic constants
$C_\varepsilon,c_\varepsilon>0$, and a positive integer
$N_\varepsilon$ such that
\begin{equation}\label{eq:intro-annealed}
 \sup_{\boldsymbol\nu\in\mathcal U_\varepsilon}
 \int_\Omega q_n^{\boldsymbol\nu}(\omega,\varepsilon)\,\dd P(\omega)
 \le C_\varepsilon e^{-c_\varepsilon n}
 \qquad(n\ge N_\varepsilon).
\end{equation}
Here $q_n^{\boldsymbol\nu}(\omega,\varepsilon)$ is the corresponding
supremal conditional deviation probability for the kernel $\boldsymbol\nu$.

The annealed deviation bound \eqref{eq:intro-annealed} enables us to study
H\"older regularity of $L_1$ under perturbations of the sampling
kernel. We use two Wasserstein distances. When the torus translation
marginals are fixed, the \emph{synchronized Wasserstein distance}
$\Wsyn$ couples cocycle elements with identical translations and
compares their matrix fields. When the translation marginals vary,
we use the ordinary Wasserstein distance $\W$ instead and require a common
phase H\"older exponent and bound for the reference and nearby matrix fields.
Assume the reference kernel $\boldsymbol\nu_0$ has an ergodic mixed base,
satisfies \SVU, and has a simple top Lyapunov exponent. We prove in
\cref{thmholder} the following local H\"older estimates, each for pairs of
kernels in a relative neighborhood of $\boldsymbol\nu_0$ within the
corresponding admissible class, with its own positive constants:
\begin{equation}\label{eq:intro-holder}
\begin{aligned}
 |L_1(\boldsymbol\nu)-L_1(\widetilde{\boldsymbol\nu})|
 &\le C\Wsyn(\boldsymbol\nu,\widetilde{\boldsymbol\nu})^\gamma,\\
 |L_1(\boldsymbol\nu)-L_1(\widetilde{\boldsymbol\nu})|
 &\le C'\W(\boldsymbol\nu,\widetilde{\boldsymbol\nu})^{\gamma'}.
\end{aligned}
\end{equation}

Duarte and Klein developed a method for H\"older regularity using
large-deviation estimates and the Avalanche Principle, with separate
applications to random and quasiperiodic cocycles \cite{DKbook}. We apply
this method to mixed random--quasiperiodic cocycles with
environment-dependent one-step laws. Related H\"older results assume uniform ergodicity of Markov
cocycles \cite{CDKM} or uniform base and fiber deviations
\cite{CaiDeng}. In contrast, our finite-state moving-law setting requires neither
uniform ergodicity of the enlarged phase--matrix dynamics nor uniform
base-and-fiber deviation estimates as additional hypotheses. Tall studies
perturbations of Markov transition probabilities with cocycle maps,
translations, and support graph fixed \cite{Tall2026}, whereas we fix
the Markov transition matrix and vary the state-indexed one-step laws.

We apply the fixed-translation Wasserstein estimate to moving-law
Schr\"odinger cocycles. For a continuous background potential $v$,
the transfer matrix at energy $E$ and noise value $\zeta$ is
\[
 A_{E,\zeta}(\theta)=
 \begin{pmatrix}v(\theta)+\zeta-E&-1\\1&0\end{pmatrix}.
\]
At each environmental state, we consider a law of frequency--noise
pairs $(\alpha,\zeta)$, whose components may be correlated.
Under suitable frequency-support and noise-nondegeneracy conditions,
we establish H\"older continuity of $E\mapsto L_1(E)$ on each compact
energy interval in \cref{thmschrodinger}.

Section~\ref{secsetup} sets up the moving-law cocycle and stationary projective
families. Section~\ref{secqualitative} develops \SVU\ and conditional expected
growth. Section~\ref{secquenched} establishes quenched concentration and the
renewal obstruction. Section~\ref{secfinite-state} establishes annealed
exponential deviations uniform over nearby sampling kernels and common
quenched starting scales for finite-state Markov environments. Section~\ref{secholder} proves the two H\"older regularity estimates,
and Section~\ref{secschrodinger} gives the Schr\"odinger application.

\section{Moving-law mixed cocycles and stationary projective families}
\label{secsetup}
We introduce the moving-law mixed cocycle model and fix the notation
for conditional growth. The stationary-family framework and basic
results from \cite{LCZpositivity} are recalled with proofs in the
present notation for later use.

\subsection{The moving-law cocycle model}
Let $(\Omega,\mathcal F,P)$ be a standard Borel probability space and
let $\Theta:\Omega\to\Omega$ be an invertible bimeasurable ergodic
probability-preserving transformation.  Fix $d\geq1$ and $r\geq2$, and
set
\[
 \cG=\TT^d\times C^0(\TT^d,\SL_r(\RR)),
\]
where the second factor has the uniform topology and $\cG$ the product
topology.  Thus $\cG$ is a Polish space.  For
$g=(\alpha,B)\in\cG$, write
\[
 a(g):=\alpha\in\TT^d,
 \qquad
 A(g):=B\in C^0(\TT^d,\SL_r(\RR)).
\]
The action on the phase--vector bundle is
\[
 g\cdot(\theta,v)
 :=\bigl(\theta+a(g),A(g)(\theta)v\bigr),
 \qquad (\theta,v)\in\TT^d\times\RR^r.
\]

Let $\nu^\bullet=(\nu_\omega)_{\omega\in\Omega}$ be a measurable
probability kernel from $\Omega$ to $\cG$. On $\cG^{\mathbb N}$,
let $\xi=(\xi_j)_{j\geq0}$ be the coordinate sequence and $\sigma$
the left shift. Conditional independence with one-step laws
$\nu_{\Theta^j\omega}$ gives the path law
\begin{equation}\label{eq:path-kernel}
 \PP_\omega=\bigotimes_{j\geq0}\nu_{\Theta^j\omega},
 \qquad
 \PP_\omega=\nu_\omega\otimes\PP_{\Theta\omega}
 \quad\text{on }\cG\times\cG^{\mathbb N}.
\end{equation}
The product-kernel construction (see \cite[Theorem~8.24]{Kallenberg2021})
makes $\omega\mapsto\PP_\omega$ a measurable probability kernel.  We write $\EE_\omega$ for integration
with respect to $\PP_\omega$.

Let $m$ be normalized Haar probability on $\TT^d$.  The mixed-base map
and its natural invariant probability are
\begin{equation}\label{eqmixed-base}
 \begin{split}
 F&:\Omega\times\cG^{\mathbb N}\times\TT^d
       \longrightarrow\Omega\times\cG^{\mathbb N}\times\TT^d,\\
 F(\omega,\xi,\theta)&=
       (\Theta\omega,\sigma\xi,\theta+a(\xi_0)),\\
 \dd\rho(\omega,\xi,\theta)&=
       \dd P(\omega)\dd\PP_\omega(\xi)\dd m(\theta).
 \end{split}
\end{equation}
The decomposition in \eqref{eq:path-kernel}, the identity
$\Theta_*P=P$, and translation invariance of $m$ give $F_*\rho=\rho$.
The \emph{mixed base} \eqref{eqmixed-base} is called ergodic when
$F$ is ergodic with respect to $\rho$.  This is a
property of the environment--path--phase system, not merely of the
environment map $\Theta$.

For $\xi\in\cG^{\mathbb N}$, define
$\Gamma_0(\xi)=0$ and
$\Gamma_n(\xi)=\sum_{j=0}^{n-1}a(\xi_j)$ for $n\geq1$.

\begin{lemma}[Dependence of mixed-base ergodicity on translations]
\label[lemma]{lemtranslation-ergodicity}
Fix the environment system $(\Omega,\mathcal F,P,\Theta)$.  If two
measurable kernels on $\cG$ have the same translation marginals
$\mu_\omega=a_*\nu_\omega$ for $P$-almost every $\omega$, then their
mixed bases are either both ergodic or both nonergodic.  In particular,
mixed-base ergodicity is common to every class with fixed translation
marginals.  If $\mu_\omega=m$ almost everywhere, then the mixed base is
ergodic.
\end{lemma}
We prove this translation-factor consequence of
\cite[Theorem~2.3]{LCZpositivity} by identifying the relevant
harmonic functions.

\begin{proof}
The translation marginals determine the Markov operator
\[
 (\mathsf Ku)(\omega,\theta)
 :=\int_{\TT^d}u(\Theta\omega,\theta+\alpha)
                         \dd\mu_\omega(\alpha)
\]
on $L^\infty(\Omega\times\TT^d,P\otimes m)$, for which $P\otimes m$
is stationary.  We first show that the mixed base is ergodic if and
only if every bounded real $\mathsf K$-harmonic function is constant.

Suppose $u=\mathsf Ku$.  Jensen's inequality and stationarity give
\[
 0\leq\int\bigl(\mathsf K(u^2)-(\mathsf Ku)^2\bigr)\dd(P\otimes m)=0.
\]
Thus equality holds in the conditional Jensen inequality, and hence
\[
 u(\Theta\omega,\theta+\alpha)=u(\omega,\theta)
\]
for $P(\dd\omega)m(\dd\theta)\mu_\omega(\dd\alpha)$-almost every
triple.  Pulling this relation back by $\alpha=a(g)$ shows that the
lift $(\omega,\xi,\theta)\mapsto u(\omega,\theta)$ is $F$-invariant.
Mixed-base ergodicity therefore forces $u$ to be constant.

Conversely, let $U$ be a bounded real $F$-invariant function and set
\[
 u(\omega,\theta)
 :=\int_{\cG^{\mathbb N}}U(\omega,\xi,\theta)\dd\PP_\omega(\xi).
\]
The product decomposition in \eqref{eq:path-kernel} and the invariance
of $U$ give $u=\mathsf Ku$.  If every bounded harmonic function is
constant, write $u=c$.  For
$\mathcal A_n=\sigma(\omega,\theta,\xi_0,\ldots,\xi_{n-1})$,
conditional independence of the tail gives
\[
 \EE_\rho[U\mid\mathcal A_n]
 =u(\Theta^n\omega,\theta+\Gamma_n(\xi))=c.
\]
Since $(\mathcal A_n)$ increases to the full sigma-algebra, upward
martingale convergence yields $U=c$ almost everywhere.  This proves
the claimed characterization, which depends only on
$\mu^\bullet$.

If $\mu_\omega=m$, then $\mathsf Ku$ is independent of $\theta$.
A harmonic $u$ therefore has the form $u(\omega,\theta)=b(\omega)$,
and the harmonic equation becomes $b=b\circ\Theta$.  Environment
ergodicity makes $b$ constant.
\end{proof}

The Lyapunov theory uses the logarithmic-integrability assumption
\begin{equation}\label{eq:log-integrability}
 \int_\Omega\int_\cG\Lambda_A(g)\dd\nu_\omega(g)\dd P(\omega)<\infty,
 \qquad
 \Lambda_A(g):=\sup_{\theta\in\TT^d}\log^+\|A(g)(\theta)\|.
\end{equation}
All matrix norms are Euclidean operator norms.  Since
$A(g)(\theta)\in\SL_r(\RR)$, its singular values give
\[
 \log^+\|A(g)(\theta)^{-1}\|
 \leq(r-1)\Lambda_A(g).
\]
Thus \eqref{eq:log-integrability} supplies the forward and inverse
integrability required by the multiplicative ergodic theorem.
 When the mixed base is ergodic, the Lyapunov exponents relative to
$\rho$ are constants
\[
 L_1(\nu^\bullet)\geq\cdots\geq L_r(\nu^\bullet),
\]
and we suppress the kernel argument when it is fixed. 

The associated matrix products are
\begin{equation}\label{eqmatrix-product} 
A_0(\xi,\theta)=I_r, 
\qquad 
A_n(\xi,\theta) 
=A(\xi_{n-1})(\theta+\Gamma_{n-1}(\xi))\cdots  
A(\xi_0)(\theta).
\end{equation}
They satisfy
\[ 
A_{n+k}(\xi,\theta) 
=A_k(\sigma^n\xi,\theta+\Gamma_n(\xi))A_n(\xi,\theta).
\]  Kingman's and
Oseledets' theorems \cite{Kingman,Oseledets} give
\begin{equation}\label{eqlyapunov}
 L_1
 =\lim_{n\to\infty}\frac1n\log\|A_n(\xi,\theta)\|
 =\lim_{n\to\infty}\frac1n
   \int\log\|A_n(\xi,\theta)\|\dd\rho(\omega,\xi,\theta),
\end{equation}
where the first convergence holds $\rho$-almost surely and in
$L^1(\rho)$.  Since the matrices have determinant one,
$\sum_{i=1}^rL_i=0$ and hence $L_1\geq0$.  For $1\leq k\leq r$,
$\wedge^k A_n$ denotes the induced linear map on
$\wedge^k\RR^r$, equipped with its Euclidean exterior-power norm.  In
particular,
\begin{equation}\label{eqexterior-exponent}
 L_1+L_2
 =\lim_{n\to\infty}\frac1n
   \int\log\|\wedge^2A_n(\xi,\theta)\|
   \dd\rho(\omega,\xi,\theta).
\end{equation}
\subsection{Projective dynamics and expected growth}
Let $Q=\TT^d\times\PP(\RR^r)$, where $\PP(\RR^r)$ is real projective
space.  Write $p\in\hat p$ for any nonzero representative.  The induced
projective action is
\[
 g\cdot(\theta,\hat p)
 :=\bigl(\theta+a(g),\widehat{A(g)(\theta)p}\bigr).
\]
The one-step growth observable
\[
 \Psi(g,(\theta,\hat p))
 :=\log\frac{\|A(g)(\theta)p\|}{\|p\|}
\]
and the projective action are well defined and continuous.

For a common full-measure domain of conditional expectations, set
\[
 h(\omega):=(r-1)\int_\cG\Lambda_A(g)\dd\nu_\omega(g).
\]
Then $h\ge0$ is measurable and integrable by
\eqref{eq:log-integrability}. Define
\[
 \Omega_*:=\bigcap_{j\in\ZZ}\Theta^{-j}\{h<\infty\}.
\]
Because $\Theta$ is invertible and preserves $P$, the set $\Omega_*$
is Borel, $\Theta$-invariant, and of full $P$-measure.

Fix $\omega\in\Omega_*$, a path $\xi=(\xi_j)_{j\ge0}$ under
$\PP_\omega$, and $x=(\theta,\hat p)\in Q$. Set $X_0^x=x$ and
$X_{j+1}^x=\xi_j\cdot X_j^x$. For $n\ge1$, define the logarithmic
growth, its conditional expectation, and, when the mixed base is
ergodic, the uniform expected-growth error by
\begin{equation}\label{eqadditive}
\begin{split}
Z_n^\omega(x)(\xi)
&:=\sum_{j=0}^{n-1}\Psi(\xi_j,X_j^x)
 =\log\frac{\|A_n(\xi,\theta)p\|}{\|p\|},\\
H_n(\omega,x)
&:=\EE_\omega Z_n^\omega(x),\\
\Delta_n(\omega)
&:=\sup_{x\in Q}
\left|\frac{H_n(\omega,x)}{n}-L_1\right|.
\end{split}
\end{equation}
The expression for $Z_n^\omega(x)$ is pathwise, with the superscript
$\omega$ indicating the conditional law $\PP_\omega$.
The path argument is usually suppressed. We extend $H_n(\omega,\cdot)$
by zero to $\Omega\setminus\Omega_*$ and, whenever $\Delta_n$ is
used, define it on all of $\Omega$ by the last line of
\eqref{eqadditive} with this extension.

The forward and inverse matrix bounds give
$|\Psi(g,x)|\leq(r-1)\Lambda_A(g)$. Hence, for
$\omega\in\Omega_*$, $n\geq1$, and every $x\in Q$,
\begin{equation}\label{eqlog-envelope}
|Z_n^\omega(x)|
\leq(r-1)\sum_{j=0}^{n-1}\Lambda_A(\xi_j),
\qquad
\|H_n(\omega,\cdot)\|_\infty
\leq\sum_{j=0}^{n-1}h(\Theta^j\omega).
\end{equation}
The second bound follows by conditional expectation, since $\xi_j$
has law $\nu_{\Theta^j\omega}$. Pathwise continuity of $Z_n(\cdot)$ and the envelope in
\eqref{eqlog-envelope} give $H_n(\omega,\cdot)\in C(Q)$
by dominated convergence. Kernel
integration gives scalar measurability in $\omega$, and the
Carath\'eodory theorem gives joint measurability. The zero extension
preserves these properties.

Separability of $C(Q)$ makes $\omega\mapsto H_n(\omega,\cdot)$
strongly measurable. Moreover,
\[
 \int_\Omega\|H_n(\omega,\cdot)\|_\infty\dd P(\omega)
 \le n\|h\|_{L^1(P)},
\]
so $H_n\in L^1(\Omega;C(Q))$. Its maximum and minimum, and hence
$\Delta_n$, are measurable because extrema over $Q$ can be taken over
one countable dense subset. Invariance of $\Omega_*$ makes the
extension consistent with all iterates.

\subsection{Stationary projective families and the Furstenberg formula}
A \emph{measurable family} of probability measures on $Q$ is a
probability kernel
$\eta^\bullet=(\eta_\omega)_{\omega\in\Omega}$.  Since $Q$ is compact
metric, this is equivalently measurability of
$\omega\mapsto\int_Q f\dd\eta_\omega$ for every $f\in C(Q)$.  Such a
family is a \emph{stationary projective family} if
\begin{equation}\label{eq:stationarity-family}
 \eta_{\Theta\omega}
 =\int_\cG g_*\eta_\omega\dd\nu_\omega(g)
 \quad\text{for }P\text{-almost every }\omega.
\end{equation}
Here $g_*\eta_\omega$ is the projective pushforward. Kernel integration
makes the right-hand side measurable, and a countable determining
family in $C(Q)$ gives one full-measure set on which the identity
holds as an equality of measures. Let
$\cS_Q(\nu^\bullet)$ denote the set of stationary families modulo
$P$-almost-everywhere equality.

Every measurable family defines the Young measure
\[
 \lambda^\eta(\dd\omega,\dd x)
 :=P(\dd\omega)\eta_\omega(\dd x)
\]
in the fixed-marginal domain
\begin{equation}\label{eq:young-domain}
 \cY_P(Q)=\{\lambda\in\Prob(\Omega\times Q):
                         (\pi_\Omega)_*\lambda=P\}.
\end{equation}
Conversely, standard-Borel disintegration \cite[Theorem~8.5]{Kallenberg2021}
represents each
$\lambda\in\cY_P(Q)$ by a measurable family, uniquely up to
$P$-almost-everywhere equality.  The stable topology
\cite{CastaingRaynaudValadier2004} is the coarsest
topology for which
$\lambda\mapsto\int\varphi\dd\lambda$ is continuous for every
$\varphi\in L^1(\Omega;C(Q))$.  Since $Q$ is compact metric and the
base is standard Borel, $\cY_P(Q)$ is compact and metrizable in this
topology \cite[Proposition~B.2]{BN2026}.

The lifted Markov operator on $L^1(\Omega;C(Q))$ is
defined by
\begin{equation}\label{eqlifted-operator}
(\widehat{\mathsf P}\varphi)(\omega,x)
=\int_\cG\varphi(\Theta\omega,g\cdot x)\dd\nu_\omega(g).
\end{equation}
Kernel integration and dominated convergence give scalar measurability
and fiberwise continuity. Since $C(Q)$ is separable, the resulting map
is strongly measurable. The bound
\[
 \|(\widehat{\mathsf P}\varphi)(\omega,\cdot)\|_\infty
 \le\|\varphi(\Theta\omega,\cdot)\|_\infty
\]
and $\Theta_*P=P$ show that $\widehat{\mathsf P}$ is a contraction
on $L^1(\Omega;C(Q))$.

The dual affine operator
$\cP:\cY_P(Q)\to\cY_P(Q)$ is defined by
\[
\int\varphi\dd(\cP\lambda)
=\int\widehat{\mathsf P}\varphi\dd\lambda,
\qquad \varphi\in L^1(\Omega;C(Q)).
\]
The transition $(\omega,x)\mapsto(\Theta\omega,g\cdot x)$ preserves
the marginal $P$, and $\cP$ is stable continuous because
$\widehat{\mathsf P}$ preserves the test space. Testing
$\cP\lambda^\eta=\lambda^\eta$ against $\chi(\omega)f(x)$ with
$\chi\in L^\infty(P)$ and $f\in C(Q)$, changing variables by $\Theta$,
and using a countable determining family
shows that its fixed points are exactly the stationary families.
Conversely, \eqref{eq:stationarity-family} and Fubini give
\[
 \int\widehat{\mathsf P}\varphi\dd\lambda^\eta
 =\int_\Omega\int_Q\varphi(\Theta\omega,x)
       \dd\eta_{\Theta\omega}(x)\dd P(\omega)
 =\int\varphi\dd\lambda^\eta.
\]

Under \eqref{eq:log-integrability}, define
\[
 \psi_{\nu^\bullet}(\omega,x)
 :=\int_\cG\Psi(g,x)\dd\nu_\omega(g),
 \qquad \omega\in\Omega_*,
\]
with the same zero extension outside $\Omega_*$. Thus
$\psi_{\nu^\bullet}=H_1\in L^1(\Omega;C(Q))$. For a stationary
family, define the Furstenberg functional
\[
 \alpha_{\nu^\bullet}(\eta^\bullet)
 :=\int_{\Omega\times Q}\psi_{\nu^\bullet}\dd\lambda^\eta.
\]
It is finite, affine, and stable continuous. 
Successive conditioning and the Markov property give
\begin{equation}\label{eqH-additive}
\begin{split}
H_n(\omega,\cdot)
&=\sum_{j=0}^{n-1}
\widehat{\mathsf P}^{\,j}\psi_{\nu^\bullet}(\omega,\cdot),\\
\int H_n\dd\lambda^\eta
&=n\alpha_{\nu^\bullet}(\eta^\bullet),
\qquad \eta^\bullet\in\cS_Q(\nu^\bullet).
\end{split}
\end{equation}
These identities use the chosen representatives and the fixed-point
relation $\cP\lambda^\eta=\lambda^\eta$.

\begin{theorem}[Stationary families and the Furstenberg formula]
\label[theorem]{thmstationary}
For the measurable kernel $\nu^\bullet$ above, the following statements
hold.
\begin{enumerate}
\item The set $\cS_Q(\nu^\bullet)$ is nonempty, convex, and stable
compact.
\item If the mixed base is ergodic, then every
$\eta^\bullet\in\cS_Q(\nu^\bullet)$ has Haar phase marginal:
$(\pi_{\TT^d})_*\eta_\omega=m$ for $P$-almost every $\omega$.
\item If the mixed base is ergodic and \eqref{eq:log-integrability}
holds, then
\[
 L_1(\nu^\bullet)
 =\max_{\eta^\bullet\in\cS_Q(\nu^\bullet)}
      \alpha_{\nu^\bullet}(\eta^\bullet).
\]
\end{enumerate}
\end{theorem}
The existence assertion in (i), and assertions (ii)--(iii), correspond
to \cite[Proposition~3.4 and Theorems~3.5 and~3.12]{LCZpositivity},
respectively.  Convexity and stable compactness in (i) follow from the
same fixed-point realization.  We include the argument in the present
notation.

\begin{proof}
For $\lambda\in\cY_P(Q)$, let
$\overline\lambda_n=n^{-1}\sum_{j=0}^{n-1}\cP^j\lambda$.  Compactness
of $\cY_P(Q)$ gives stable cluster points.  For every
$\varphi\in L^1(\Omega;C(Q))$,
\[
 \left|\int\varphi\dd(\cP\overline\lambda_n-
                              \overline\lambda_n)\right|
 =\frac1n\left|\int\varphi\dd(\cP^n\lambda-\lambda)\right|
 \leq\frac{2}{n}\|\varphi\|_{L^1(\Omega;C(Q))}.
\]
Every cluster point is therefore fixed by $\cP$.  The fixed-point set
is nonempty, closed, and convex in the stable compact space
$\cY_P(Q)$, proving (i).

For (ii), let
$\lambda_\omega=(\pi_{\TT^d})_*\eta_\omega$ and
$\mu_\omega=a_*\nu_\omega$.  Projecting
\eqref{eq:stationarity-family} onto the torus gives
$\lambda_{\Theta\omega}=\mu_\omega*\lambda_\omega$.  For
$k\in\ZZ^d\setminus\{0\}$, put
\[
 c_k(\omega)=\int_{\TT^d}e^{2\pi i k\cdot\theta}
                         \dd\lambda_\omega(\theta),
 \qquad
 b_k(\omega)=\int_{\TT^d}e^{2\pi i k\cdot\alpha}
                         \dd\mu_\omega(\alpha).
\]
Then $c_k\circ\Theta=b_k c_k$ and
$|c_k\circ\Theta|\leq|c_k|$.  Since $\Theta$ preserves $P$,
integration forces $|c_k\circ\Theta|=|c_k|$ almost everywhere.
Thus $|c_k|$ is $\Theta$-invariant and hence constant by environment
ergodicity.  If the constant were positive, then
$|b_k|=1$ almost everywhere, and equality in the triangle inequality
would give
$e^{2\pi i k\cdot a(g)}=b_k(\omega)$ for
$\nu_\omega$-almost every $g$.  The function
\[
 (\omega,\xi,\theta)\longmapsto
 \frac{e^{2\pi i k\cdot\theta}}{c_k(\omega)}
\]
would then be a nonconstant $F$-invariant function, contradicting
mixed-base ergodicity.  Hence all nonzero Fourier coefficients of
$\lambda_\omega$ vanish, and Fourier uniqueness gives
$\lambda_\omega=m$ almost everywhere.

For (iii), let $m_{\PP}$ be the pushforward to $\PP(\RR^r)$ of
normalized rotation-invariant probability on the Euclidean unit
sphere.  There is $C_r<\infty$ such that every invertible matrix $T$
satisfies
\begin{equation}\label{eqprojective-comparison}
 \log\|T\|-C_r
 \leq\int_{\PP(\RR^r)}
       \log\frac{\|Tp\|}{\|p\|}\dd m_{\PP}(\hat p)
 \leq\log\|T\|.
\end{equation}
Indeed, if $u$ is a unit top right singular vector, then
$\|Tp\|\geq\|T\||\langle p,u\rangle|$ for unit $p$, and the integral
of $\log|\langle p,u\rangle|$ is finite and independent of $u$.

Let $\lambda_0=P\otimes m\otimes m_{\PP}\in\cY_P(Q)$ and set
\[
 \overline\lambda_n^{\,0}
 :=\frac1n\sum_{j=0}^{n-1}\cP^j\lambda_0.
\]
Choose a subsequence $n_j$ along which
$\overline\lambda_{n_j}^{\,0}\to\lambda^\eta$ stably.  By (i),
$\eta^\bullet\in\cS_Q(\nu^\bullet)$.  Equations
\eqref{eqH-additive}, \eqref{eqprojective-comparison}, and
\eqref{eqlyapunov} give
\[
 \int\psi_{\nu^\bullet}\dd\overline\lambda_{n_j}^{\,0}
 =\frac1{n_j}\int H_{n_j}\dd\lambda_0
 \longrightarrow L_1.
\]
Stable convergence against
$\psi_{\nu^\bullet}\in L^1(\Omega;C(Q))$ yields
$\alpha_{\nu^\bullet}(\eta^\bullet)=L_1$.

Conversely, if $\eta^\bullet$ is stationary, then (ii) and
\eqref{eqH-additive} imply
\[
 n\alpha_{\nu^\bullet}(\eta^\bullet)
 =\int\log\frac{\|A_n(\xi,\theta)p\|}{\|p\|}
       \dd P(\omega)\dd\PP_\omega(\xi)
       \dd\eta_\omega(\theta,\hat p)
 \leq\int\log\|A_n(\xi,\theta)\|\dd\rho.
\]
Divide by $n$ and use \eqref{eqlyapunov}.  This gives the reverse bound
for every stationary family and proves the variational formula.
\end{proof}

\begin{definition}[Stationary-value uniqueness]
\label[definition]{defSVU}
Assume mixed-base ergodicity and \eqref{eq:log-integrability}.
We say that the kernel $\nu^\bullet$ has
\emph{stationary-value uniqueness} (\SVU) if
\[
\alpha_{\nu^\bullet}(\eta^\bullet)=L_1(\nu^\bullet)
\qquad\text{for every }
\eta^\bullet\in\cS_Q(\nu^\bullet).
\]
\end{definition}

We record an integral-envelope identity for additive Markov sums.

\begin{lemma}[Integral envelopes of additive Markov sums]
\label[lemma]{lem:integrated-envelopes}
Let $\varphi\in L^1(\Omega;C(Q))$ be real-valued and set
\[
 S_n\varphi=\sum_{j=0}^{n-1}\widehat{\mathsf P}^{\,j}\varphi.
\]
Then
\[
\begin{aligned}
 \lim_{n\to\infty}\frac1n
 \int_\Omega\max_{x\in Q}(S_n\varphi)(\omega,x)\dd P(\omega)
 &=\max_{\eta^\bullet\in\cS_Q(\nu^\bullet)}
       \int_{\Omega\times Q}\varphi\dd\lambda^\eta,\\
 \lim_{n\to\infty}\frac1n
 \int_\Omega\min_{x\in Q}(S_n\varphi)(\omega,x)\dd P(\omega)
 &=\min_{\eta^\bullet\in\cS_Q(\nu^\bullet)}
       \int_{\Omega\times Q}\varphi\dd\lambda^\eta.
\end{aligned}
\]
Both extrema over stationary families are attained.
\end{lemma}

\begin{proof}
Stable compactness of $\cS_Q(\nu^\bullet)$ and the definition of the
stable topology give attainment of the extrema. It suffices to prove
the maximum identity, since the minimum identity follows by replacing
$\varphi$ with $-\varphi$.

Write
\[
 M_n(\omega)=\max_{x\in Q}(S_n\varphi)(\omega,x),
 \qquad
 \beta=\max_{\eta^\bullet\in\cS_Q(\nu^\bullet)}
                 \int\varphi\dd\lambda^\eta.
\]
Since $\widehat{\mathsf P}$ preserves $L^1(\Omega;C(Q))$ and is a
contraction there, $M_n$ is measurable and
\[
 \|M_n\|_{L^1(P)}
 \leq\|S_n\varphi\|_{L^1(\Omega;C(Q))}
 \leq n\|\varphi\|_{L^1(\Omega;C(Q))}.
\]
For every stationary family, duality and
$\cP\lambda^\eta=\lambda^\eta$ give
\[
 n\int\varphi\dd\lambda^\eta
 =\int S_n\varphi\dd\lambda^\eta
 \leq\int_\Omega M_n\dd P.
\]
Thus $\beta\leq n^{-1}\int M_n\dd P$ for every $n$.

For the reverse limiting inequality, fix a countable dense set
$(q_k)_{k\geq1}$ in $Q$. On a common full-measure set where all
$S_n\varphi$ have continuous fibers, define
\[
 k_n(\omega)=\min\bigl\{k\geq1:
 (S_n\varphi)(\omega,q_k)>M_n(\omega)-1\bigr\},
 \qquad x_n(\omega)=q_{k_n(\omega)}.
\]
Continuity and density make $k_n$ finite, and the defining measurable
inequalities make $k_n$ measurable. Extend $x_n$ by $q_1$ on the
exceptional set. Define
\[
 \lambda_n(\dd\omega,\dd x)
 =P(\dd\omega)\delta_{x_n(\omega)}(\dd x),
 \qquad
 \overline\lambda_n=\frac1n\sum_{j=0}^{n-1}\cP^j\lambda_n.
\]
Both measures belong to $\cY_P(Q)$. For every
$\chi\in L^1(\Omega;C(Q))$, telescoping and the common environment
marginal give
\[
\begin{aligned}
 \left|\int\chi\dd(\cP\overline\lambda_n-\overline\lambda_n)\right|
 &=\frac1n\left|\int\chi\dd(\cP^n\lambda_n-\lambda_n)\right|\\
 &\leq\frac2n\|\chi\|_{L^1(\Omega;C(Q))}.
\end{aligned}
\]
Hence every stable subsequential limit of $\overline\lambda_n$ is
fixed by $\cP$, by stable continuity of $\cP$, and therefore equals
$\lambda^\eta$ for some $\eta^\bullet\in\cS_Q(\nu^\bullet)$.
Moreover, duality and the choice of $x_n$ imply
\[
 0\leq\frac1n\int_\Omega M_n\dd P
       -\int\varphi\dd\overline\lambda_n
 =\frac1n\int_\Omega
   \bigl[M_n(\omega)-(S_n\varphi)(\omega,x_n(\omega))\bigr]
   \dd P(\omega)
 \leq\frac1n.
\]
Choose a subsequence along which $n^{-1}\int M_n\dd P$ converges to
its limsup. Stable compactness of $\cY_P(Q)$ gives a further
subsequence for which $\overline\lambda_n\to\lambda^\eta$ stably.
Since $\varphi$ belongs to the stable test space, the preceding
estimate yields
\[
 \limsup_{n\to\infty}\frac1n\int_\Omega M_n\dd P
 =\int\varphi\dd\lambda^\eta\leq\beta.
\]
Together with the opposite inequality, this proves the claim.
\end{proof}

\section{Stationary values and qualitative growth}\label{secqualitative}
Under mixed-base ergodicity and logarithmic integrability, relative
quasi-irreducibility is equivalent to SVU and forces uniform conditional expected growth.
This yields qualitative continuity at a reference SVU kernel and the
small block drift needed for concentration. For compactly supported
fixed laws, the stationary-value criterion and a uniform expected-growth
theorem with a top Lyapunov gap appear in
\cite[Proposition~3.7 and Theorem~1.2]{CDKfurstenberg}.

\subsection{Relative quasi-irreducibility and SVU}
A submaximal stationary value produces a relative invariant subbundle
with submaximal growth, and conversely such a subbundle supports a
stationary family with submaximal value. This correspondence proves
\cref{thm:QI-SVU-equivalence}. 

\begin{definition}[Relative invariant subbundle]\label[definition]{def:relative-invariant-subbundle}
Let $1\leq k<r$.  A measurable map
\(
 L:\Omega\times\TT^d\longrightarrow\Gr_k(\RR^r)
\)
is a $k$-dimensional relative invariant subbundle if
\begin{equation}\label{eq:a.e.-bundle-invariance}
 A(g)(\theta)L(\omega,\theta)
 =L(\Theta\omega,\theta+a(g))
\end{equation}
for $P(\dd\omega)\nu_\omega(\dd g)m(\dd\theta)$-almost every $(\omega,g,\theta)$.
\end{definition}

\begin{remark}[Measurable subbundle representations]
\label[remark]{rem:equivalent-subbundle-descriptions}
Relative invariant subbundles have fixed rank $1\le k<r$.
A \emph{measurable linear-subspace field} may instead have rank $0$,
rank $r$, or varying rank. Appendix~\ref{app:subbundle-coding}
identifies Grassmannian maps with measurable total spaces and supplies
measurable frames and fixed-rank versions.
\end{remark}

For such a bundle, define
\begin{equation}\label{eq:restricted-subadditive}
 f_n^L(\omega,\xi,\theta)
 :=\log\left\|A_n(\xi,\theta)|_{L(\omega,\theta)}\right\|.
\end{equation}

\begin{proposition}[Restricted top Lyapunov exponent]\label[proposition]{prop:restricted-exponent}
Assume mixed-base ergodicity and \eqref{eq:log-integrability}.  If $L$ is a relative invariant subbundle, then there is a constant $\lambda_1(L)\in\RR$ such that
\begin{equation}\label{eq:restricted-limit}
 \lim_{n\to\infty}\frac1n f_n^L(\omega,\xi,\theta)
 =\lambda_1(L)
\end{equation}
for $\rho$-almost every $(\omega,\xi,\theta)$.  Moreover,
\begin{equation}\label{eq:restricted-integral-limit}
 \lambda_1(L)
 =\lim_{n\to\infty}\frac1n\int f_n^L\dd\rho
 =\inf_{n\geq1}\frac1n\int f_n^L\dd\rho.
\end{equation}
\end{proposition}

\begin{proof}
The map $(B,W)\mapsto\|B|_W\|$ is continuous on $\GL_r(\RR)\times\Gr_k(\RR^r)$, so $f_n^L$ is measurable.  By Fubini and the conditional product structure \eqref{eq:path-kernel}, after discarding a $\rho$-null set, the one-step invariance \eqref{eq:a.e.-bundle-invariance} holds at every time along the sampled path.  Consequently,
\(
 A_n(\xi,\theta)L(\omega,\theta)
 =L(\Theta^n\omega,\theta+\Gamma_n(\xi))
\)
for all $n$ on that full-measure set.  The cocycle identity implies
\(
 f_{n+m}^L
 \leq f_n^L+f_m^L\circ F^n.
\)
The inverse estimate following \eqref{eq:log-integrability} gives
\[
 -(r-1)\sum_{j=0}^{n-1}\Lambda_A(\xi_j)
 \leq f_n^L(\omega,\xi,\theta)
 \leq\sum_{j=0}^{n-1}\Lambda_A(\xi_j).
\]
Thus $f_n^L\in L^1(\rho)$ and
$\inf_{n\geq1}n^{-1}\int f_n^L\dd\rho> -\infty$.
Kingman's subadditive ergodic theorem \cite[Theorem~25.20]{Kallenberg2021} gives
\eqref{eq:restricted-limit} and \eqref{eq:restricted-integral-limit},
with a finite limit that is constant by mixed-base ergodicity.
\end{proof}

\begin{definition}[Relative irreducibility and quasi-irreducibility]\label[definition]{def:moving-quasi-irreducibility}
The moving cocycle is \emph{relatively irreducible} if it admits no
relative invariant subbundle.  It is \emph{relatively
quasi-irreducible} if no relative invariant subbundle $L$ satisfies
\(
  \lambda_1(L)<L_1(\nu^\bullet).
\)
\end{definition}

To relate stationary values to slow subbundles, associate to
$\eta^\bullet\in\cS_Q(\nu^\bullet)$ the measure
\begin{equation}\label{eq:projective-extension-measure}
 \dd\rho_\eta(\omega,\xi,x)
 :=\dd P(\omega)\dd\PP_\omega(\xi)\dd\eta_\omega(x)
\end{equation}
on $\Omega\times\cG^{\mathbb N}\times Q$, and let
\begin{equation}\label{eq:projective-skew-map}
 \widehat F_Q(\omega,\xi,x)
 :=(\Theta\omega,\sigma\xi,\xi_0\cdot x).
\end{equation}

\begin{lemma}[Ergodic extensions of extremal stationary families]
\label[lemma]{lem:extremal-ergodic-extension}
The measure $\rho_\eta$ is $\widehat F_Q$-invariant.  If $\eta^\bullet$ is an extreme point of $\cS_Q(\nu^\bullet)$, then $(\widehat F_Q,\rho_\eta)$ is ergodic.  In that case,
\begin{equation}\label{eq:extremal-directional-limit}
 \lim_{n\to\infty}
 \frac1n\log\frac{\|A_n(\xi,\theta)p\|}{\|p\|}
 =\alpha_{\nu^\bullet}(\eta^\bullet)
\end{equation}
for $\rho_\eta$-almost every $(\omega,\xi,\theta,\hat p)$.
\end{lemma}

\begin{proof}
For bounded measurable $\Phi$, the product-kernel identity and
\eqref{eq:stationarity-family} give
\[
 \int\Phi\circ\widehat F_Q\dd\rho_\eta
 =\int_\Omega\int_Q\int_{\cG^{\mathbb N}}
 \Phi(\Theta\omega,\zeta,y)
 \dd\PP_{\Theta\omega}(\zeta)\dd\eta_{\Theta\omega}(y)\dd P(\omega)
 =\int\Phi\dd\rho_\eta,
\]
where the last equality uses $\Theta_*P=P$.
On $Y:=\Omega\times Q$, consider the Markov operator
\[
 (\mathsf K h)(\omega,x)
 :=\int_{\cG}h(\Theta\omega,g\cdot x)\dd\nu_\omega(g).
\]
The measure $\lambda^\eta(\dd\omega,\dd x) :=P(\dd\omega)\eta_\omega(\dd x)$ is $\mathsf K$-stationary.  It is extreme among all $\mathsf K$-stationary probability measures.  Indeed, if
\(
 \lambda^\eta=t\lambda_1+(1-t)\lambda_2,
 \ 0<t<1,
\)
with $\lambda_i$ stationary, then the environment marginals $P_i$ of $\lambda_i$ are $\Theta$-invariant and satisfy $P=tP_1+(1-t)P_2$.  Since $P$ is ergodic, it is extreme among the $\Theta$-invariant probabilities, so $P_1=P_2=P$.  The measures $\lambda_i$ therefore disintegrate into stationary projective families, and extremality of $\eta^\bullet$ forces $\lambda_1=\lambda_2=\lambda^\eta$.

Every bounded real $\mathsf K$-harmonic function is constant $\lambda^\eta$-almost everywhere.  If $h=\mathsf K h$, Jensen's inequality and stationarity give
\[
 0\leq
 \int_Y\bigl(\mathsf K(h^2)-(\mathsf K h)^2\bigr)\dd\lambda^\eta
 =0.
\]
Thus $h(y')=h(y)$ for $\lambda^\eta(\dd y)\mathsf K(y,\dd y')$-almost every $(y,y')$.  Every sublevel set $B$ of $h$ is consequently invariant modulo $\lambda^\eta$, so $\mathsf K\one_B=\one_B$ almost everywhere and
\(
 \one_B\,\mathsf K f=\mathsf K(f\one_B)\  \lambda^\eta\text{-a.e.}
\)
for every bounded measurable $f$.  If $0<\lambda^\eta(B)<1$, this identity and stationarity show that the normalized restrictions of $\lambda^\eta$ to $B$ and $Y\setminus B$ are two distinct stationary probabilities, contradicting extremality.  Hence $h$ is constant.

For a bounded real $\widehat F_Q$-invariant function $H$, set
\[
 h(\omega,x)=\int H(\omega,\xi,x)\dd\PP_\omega(\xi).
\]
The product decomposition gives $h=\mathsf Kh$, hence $h=c$ almost
everywhere. With $\mathcal A_n=\sigma(\omega,x,\xi_0,\ldots,\xi_{n-1})$
and $X_n=\xi_{n-1}\cdots\xi_0\cdot x$, conditional independence gives
\[
 \EE_{\rho_\eta}[H\mid\mathcal A_n]
 =h(\Theta^n\omega,X_n)=c,
\]
since $(\Theta^n\omega,X_n)$ has law $\lambda^\eta$.
Upward martingale convergence \cite[Theorem~9.24]{Kallenberg2021}
yields $H=c$, proving ergodicity.

Finally, the observable $(\omega,\xi,x)\mapsto\Psi(\xi_0,x)$ is integrable by \eqref{eq:log-integrability}.  Its integral against $\rho_\eta$ is $\alpha_{\nu^\bullet}(\eta^\bullet)$, and the projective cocycle identity gives
\[
 \sum_{j=0}^{n-1}\Psi\circ\widehat F_Q^j
 =\log\frac{\|A_n(\xi,\theta)p\|}{\|p\|}.
\]
Birkhoff's theorem \cite[Theorem~25.6]{Kallenberg2021} proves \eqref{eq:extremal-directional-limit}.
\end{proof}

For $c\in\RR$, $v\neq0$, and $(\omega,\theta)\in\Omega\times\TT^d$, set
\begin{equation}\label{eq:slow-probability}
 p_c(\omega,\theta,v)
 :=\PP_\omega\left\{
 \xi:
 \limsup_{n\to\infty}
 \frac1n\log\frac{\|A_n(\xi,\theta)v\|}{\|v\|}
 \leq c
 \right\},
\end{equation}
and put $p_c(\omega,\theta,0):=1$.  Define
\begin{equation}\label{eq:slow-space}
 L_c(\omega,\theta)
 :=\{v\in\RR^r:p_c(\omega,\theta,v)=1\}.
\end{equation}

\begin{lemma}[Measurable slow subspaces]\label[lemma]{lem:slow-subspaces}
Assume mixed-base ergodicity and \eqref{eq:log-integrability}.  For every $c\in\RR$, the family $L_c$ is, after changing it on a $(P\times m)$-null set, a measurable linear-subspace field of constant dimension.  It satisfies
\begin{equation}\label{eq:slow-space-invariance}
 A(g)(\theta)L_c(\omega,\theta)
 =L_c(\Theta\omega,\theta+a(g))
\end{equation}
for $P(\dd\omega)m(\dd\theta)\nu_\omega(\dd g)$-almost every $(\omega,\theta,g)$.  If its common dimension $k$ satisfies $1\leq k<r$, then $L_c$ is a relative invariant subbundle and
\begin{equation}\label{eq:slow-restricted-bound}
 \lambda_1(L_c)\leq c.
\end{equation}
If $k=r$, then $L_1(\nu^\bullet)\leq c$.
\end{lemma}

\begin{proof}
For $v\neq0$, the event in \eqref{eq:slow-probability} is
\[
 \bigcap_{\ell\geq1}\bigcup_{N\geq1}\bigcap_{n\geq N}
 \left\{
 \frac1n\log\frac{\|A_n(\xi,\theta)v\|}{\|v\|}
 \leq c+\frac1\ell
 \right\}.
\]
The matrix products are jointly measurable in $(\omega,\xi,\theta,v)$ and $\omega\mapsto\PP_\omega$ is a measurable probability kernel.  Hence $p_c$ is measurable.  If $v,w\in L_c(\omega,\theta)$, then the two slow events have probability one, and on their intersection
\[
 \|A_n(\xi,\theta)(v+w)\|
 \leq\|A_n(\xi,\theta)v\|+\|A_n(\xi,\theta)w\|.
\]
It follows that $v+w\in L_c(\omega,\theta)$, and scalar multiplication is immediate.  Thus each fiber is a linear, hence closed, subspace.  The graph
\(
 \operatorname{graph}(L_c)
 :=\{(\omega,\theta,v):p_c(\omega,\theta,v)=1\}
\)
is measurable. By Appendix~\ref{app:subbundle-coding}, $L_c$ is
therefore a measurable map into
$\bigsqcup_{\ell=0}^r\Gr_\ell(\RR^r)$ and admits measurable
orthonormal frames on each dimension stratum.

For $\xi=(g,\zeta)$ and $n\geq2$, the cocycle identity gives
\[
 A_n(\xi,\theta)v
 =A_{n-1}(\zeta,\theta+a(g))A(g)(\theta)v.
\]
The one-step multiplicative factor changes the normalized logarithm by a term tending to zero.  Using $\PP_\omega=\nu_\omega\otimes\PP_{\Theta\omega}$ therefore yields
\begin{equation}\label{eq:slow-markov-identity}
 p_c(\omega,\theta,v)
 =\int_{\cG}
 p_c\bigl(\Theta\omega,\theta+a(g),A(g)(\theta)v\bigr)
 \dd\nu_\omega(g).
\end{equation}
If $v\in L_c(\omega,\theta)$, the left-hand side equals one and the integrand belongs to $[0,1]$. Hence, the transformed vector belongs to the next slow space for $\nu_\omega$-almost every $g$.  Applying this to measurable orthonormal frames on the finitely many dimension strata and taking finite intersections gives
\begin{equation}\label{eq:slow-inclusion}
 A(g)(\theta)L_c(\omega,\theta)
 \subset L_c(\Theta\omega,\theta+a(g))
\end{equation}
for $P(\dd\omega)m(\dd\theta)\nu_\omega(\dd g)$-almost every triple.

Let $d_c(\omega,\theta):=\dim L_c(\omega,\theta)$.  Since $A(g)(\theta)$ is invertible, \eqref{eq:slow-inclusion} gives
\(
 d_c(\omega,\theta)
 \leq d_c(\Theta\omega,\theta+a(g)).
\)
The measure $P\times m$ is stationary for the environment--phase transition, so the integrals of the two sides with respect to $P(\dd\omega)m(\dd\theta)\nu_\omega(\dd g)$ are equal.  The nonnegative integer-valued difference therefore vanishes almost everywhere.  Thus \eqref{eq:slow-inclusion} is an equality.  The lift $(\omega,\xi,\theta)\mapsto d_c(\omega,\theta)$ is then $F$-invariant, so mixed-base ergodicity makes it equal to a constant $k$ almost everywhere.  Assigning a fixed $W_*\in\Gr_k(\RR^r)$ on the exceptional set
gives a pointwise fixed-rank version without changing these
almost-everywhere conclusions.

Assume $k\geq1$ and choose a measurable orthonormal frame $e_1,\ldots,e_k$ of $L_c$.  By the definition of $L_c$, for $\rho$-almost every $(\omega,\xi,\theta)$ all the finitely many vectors $e_i(\omega,\theta)$ satisfy the upper-growth bound $c$.  Since
\[
 \left\|A_n(\xi,\theta)|_{L_c(\omega,\theta)}\right\|
 \leq\sqrt{k}\max_{1\leq i\leq k}
 \|A_n(\xi,\theta)e_i(\omega,\theta)\|,
\]
we obtain
\(
 \limsup_{n\to\infty}\frac1n
 \log\left\|A_n(\xi,\theta)|_{L_c(\omega,\theta)}\right\|
 \leq c
\)
for $\rho$-almost every point.  If $k<r$, Proposition~\ref{prop:restricted-exponent} gives \eqref{eq:slow-restricted-bound}.  If $k=r$, the restriction is the full matrix product, and \eqref{eqlyapunov} gives $L_1(\nu^\bullet)\leq c$.
\end{proof}

\begin{proposition}[Submaximal stationary values and slow subbundles]\label[proposition]{propslow-obstruction}
Assume mixed-base ergodicity and \eqref{eq:log-integrability}.  If an extremal stationary family $\eta^\bullet\in\cS_Q(\nu^\bullet)$ satisfies
\(
 \alpha_{\nu^\bullet}(\eta^\bullet)<L_1(\nu^\bullet),
\)
then there exists a relative invariant subbundle $L$ such that
\(
 \lambda_1(L)<L_1(\nu^\bullet).
\)
\end{proposition}

\begin{proof}
Choose $\alpha_{\nu^\bullet}(\eta^\bullet)<c<L_1(\nu^\bullet)$ and consider $L_c$ from Lemma~\ref{lem:slow-subspaces}.  By Lemma~\ref{lem:extremal-ergodic-extension} and Fubini, for $P(\dd\omega)\eta_\omega(\dd\theta,\dd\hat p)$-almost every $(\omega,\theta,\hat p)$, the direction $\hat p$ belongs to $\PP(L_c(\omega,\theta))$.  Theorem~\ref{thmstationary}(ii) gives the Haar torus marginal.  After disintegrating $\eta_\omega$ over $m$, it follows that $L_c(\omega,\theta)\neq\{0\}$ for $(P\times m)$-almost every $(\omega,\theta)$.

By \cref{lem:slow-subspaces}, full rank would imply $L_1\le c$.
Thus the common rank $k$ satisfies $1\le k<r$, and the same lemma
gives relative invariance and
$\lambda_1(L_c)\le c<L_1$. Take $L=L_c$.
\end{proof}

For a $k$-dimensional relative invariant subbundle $L$, retain support
under stable limits by adjoining the subspace coordinate in
\begin{equation}\label{eq:incidence-space}
 \mathcal I_k
 :=\left\{
 (\theta,W,\hat p)\in\TT^d\times\Gr_k(\RR^r)\times\PP(\RR^r):
 \hat p\in\PP(W)
 \right\}.
\end{equation}
The action
\begin{equation}\label{eq:incidence-action}
 g\cdot(\theta,W,\hat p)
 :=\bigl(\theta+a(g),A(g)(\theta)W,\widehat{A(g)(\theta)p}\bigr)
\end{equation}
is continuous on $\mathcal I_k$.

\begin{lemma}[Stationary families on invariant subbundles]
\label[lemma]{lem:bundle-stationary-family}
Assume mixed-base ergodicity and \eqref{eq:log-integrability}.  If $L$ is a $k$-dimensional relative invariant subbundle, then there is a stationary family $\eta_L^\bullet\in\cS_Q(\nu^\bullet)$ such that, after disintegration over the Haar torus marginal,
\begin{equation}\label{eq:bundle-supported-family}
 \eta_{L,\omega,\theta}\bigl(\PP(L(\omega,\theta))\bigr)=1
\end{equation}
for $(P\times m)$-almost every $(\omega,\theta)$.  Moreover,
\begin{equation}\label{eq:bundle-value-bound}
 \alpha_{\nu^\bullet}(\eta_L^\bullet)
 \leq\lambda_1(L).
\end{equation}
\end{lemma}

\begin{proof}
Let
\(
 \gamma_L(\dd\omega,\dd\theta,\dd W)
 :=P(\dd\omega)m(\dd\theta)
   \delta_{L(\omega,\theta)}(\dd W).
\)
Let $\cY_P(\mathcal I_k)$ denote the stable Young-measure space defined as in \eqref{eq:young-domain}, with $Q$ replaced by the compact space $\mathcal I_k$.  Consider its subset
\(
 \mathcal Y_L
 :=\left\{
 \Lambda\in\cY_P(\mathcal I_k):
 (\pi_{\Omega,\theta,W})_*\Lambda=\gamma_L
 \right\}.
\)
This set is convex and stable compact. It is nonempty: a measurable
unit section $e_1(\omega,\theta)$ of $L$ yields the measure
\(
 P(\dd\omega)m(\dd\theta)
 \delta_{L(\omega,\theta)}(\dd W)
 \delta_{\widehat{e_1(\omega,\theta)}}(\dd\hat p)
\)
in $\mathcal Y_L$. Closedness follows because the fixed $(\omega,\theta,W)$-marginal is tested by functions in $L^1(\Omega;C(\TT^d\times\Gr_k(\RR^r)))$, viewed as functions on $\mathcal I_k$ independent of $\hat p$.

Define the affine operator $\mathcal P_L$ by
\[
 \int\phi\dd(\mathcal P_L\Lambda)
 =\int_\Omega\int_{\mathcal I_k}\int_\cG
 \phi(\Theta\omega,g\cdot y)
 \dd\nu_\omega(g)\dd\Lambda_\omega(y)\dd P(\omega).
\]
For $\phi\in L^1(\Omega;C(\mathcal I_k))$, the continuous action
\eqref{eq:incidence-action} and the bound
\[
 \left\|\int_\cG\phi(\Theta\omega,g\cdot\,\bullet)
                    \dd\nu_\omega(g)\right\|_\infty
 \le\|\phi(\Theta\omega,\cdot)\|_\infty
\]
give stable continuity exactly as in \eqref{eqlifted-operator}.

For $\chi\in L^1(\Omega;C(\TT^d\times\Gr_k(\RR^r)))$ and
$\Lambda\in\mathcal Y_L$, relative invariance gives
\begin{align*}
 &\int \chi(\omega,\theta,W)
 \dd(\pi_{\Omega,\theta,W})_*(\mathcal P_L\Lambda)\\
 &\quad=\int_\Omega\int_{\TT^d}\int_{\cG}
 \chi\bigl(\Theta\omega,\theta+a(g),A(g)(\theta)L(\omega,\theta)\bigr)
 \dd\nu_\omega(g)\dd m(\theta)\dd P(\omega)\\
 &\quad=\int_\Omega\int_{\TT^d}\int_{\cG}
 \chi\bigl(\Theta\omega,\theta+a(g),L(\Theta\omega,\theta+a(g))\bigr)
 \dd\nu_\omega(g)\dd m(\theta)\dd P(\omega)\\
 &\quad=\int_\Omega\int_{\TT^d}
 \chi\bigl(\omega,\theta,L(\omega,\theta)\bigr)
 \dd m(\theta)\dd P(\omega).
\end{align*}
The last equality uses invariance of $P\times m$ under the environment--phase transition.  Hence $(\pi_{\Omega,\theta,W})_*(\mathcal P_L\Lambda)=\gamma_L$, so $\mathcal P_L(\mathcal Y_L)\subset\mathcal Y_L$.  The Ces\`aro fixed-point argument used in \cref{thmstationary}, now inside the compact invariant set $\mathcal Y_L$, gives $\Lambda_L\in\mathcal Y_L$ with $\mathcal P_L\Lambda_L=\Lambda_L$.

Project $\Lambda_L$ onto $(\omega,\theta,\hat p)$ and disintegrate the resulting Young measure over $P$. Denote its conditional projective law by $\eta_{L,\omega}$. The fixed-point identity gives $\eta_L^\bullet\in\cS_Q(\nu^\bullet)$. The fixed marginal and incidence condition give \eqref{eq:bundle-supported-family}.

Let $\rho_{\eta_L}$ be the invariant projective extension from \eqref{eq:projective-extension-measure}.  Its $(\omega,\xi,\theta)$-marginal is $\rho$, since the torus marginal of $\eta_L^\bullet$ is $m$.  Invariance and the projective cocycle identity give, for every $n\geq1$,
\begin{align*}
 n\,\alpha_{\nu^\bullet}(\eta_L^\bullet)
 &=\int
 \log\frac{\|A_n(\xi,\theta)p\|}{\|p\|}
 \dd\rho_{\eta_L}(\omega,\xi,\theta,\hat p)\\
 &\leq\int
 \log\left\|A_n(\xi,\theta)|_{L(\omega,\theta)}\right\|
 \dd\rho(\omega,\xi,\theta)\\
 &=\int f_n^L\dd\rho.
\end{align*}
The inequality uses \eqref{eq:bundle-supported-family}.
Dividing by $n$ and taking the infimum in
\eqref{eq:restricted-integral-limit} proves
\eqref{eq:bundle-value-bound}.
\end{proof}

\begin{theorem}[Relative quasi-irreducibility and stationary-value uniqueness]
\label[theorem]{thm:QI-SVU-equivalence}
Assume mixed-base ergodicity and \eqref{eq:log-integrability}.  The following conditions are equivalent:
\begin{enumerate}[label=\textup{(\roman*)}]
\item the cocycle is relatively quasi-irreducible;
\item every extreme point $\eta^\bullet$ of $\cS_Q(\nu^\bullet)$ satisfies
\(
 \alpha_{\nu^\bullet}(\eta^\bullet)=L_1(\nu^\bullet);
\)
\item the cocycle has stationary-value uniqueness.
\end{enumerate}
\end{theorem}

\begin{proof}
Condition (iii) implies (ii). If SVU fails, compactness of
$\cS_Q(\nu^\bullet)$ and continuity of $\alpha_{\nu^\bullet}$ give
an attained minimum
\(
 a_*:=\min_{\eta^\bullet\in\cS_Q(\nu^\bullet)}
 \alpha_{\nu^\bullet}(\eta^\bullet)
\)
with $a_*<L_1(\nu^\bullet)$. The minimizer set
\(
 \mathcal F_*:=\left\{
 \eta^\bullet\in\cS_Q(\nu^\bullet):
 \alpha_{\nu^\bullet}(\eta^\bullet)=a_*
 \right\}
\)
is a nonempty compact face because $\alpha_{\nu^\bullet}$ is affine.  By the Krein--Milman theorem \cite{KreinMilman1940}, $\mathcal F_*$ has an extreme point, and every extreme point of a face is extreme in the ambient compact convex set.  Thus (ii) fails.  This proves (ii)$\Leftrightarrow$(iii).

Proposition~\ref{propslow-obstruction} gives
(i)$\Rightarrow$(ii), since a submaximal extremal value would produce a
relative invariant subbundle with restricted exponent below $L_1$.

Conversely, if (i) fails, there is a relative invariant subbundle $L$
with $\lambda_1(L)<L_1$.  Lemma~\ref{lem:bundle-stationary-family}
gives a stationary family $\eta_L^\bullet$ satisfying
\(
 \alpha_{\nu^\bullet}(\eta_L^\bullet)
 \leq\lambda_1(L)<L_1(\nu^\bullet),
\)
so stationary-value uniqueness fails.  This proves (iii)$\Rightarrow$(i).
\end{proof}

\begin{definition}[Environment-wise strong irreducibility]
A finite family of nonzero proper measurable subbundles is
\emph{support-invariant} if its fiberwise union $U(\omega,\theta)$ satisfies
\[
 A(g)(\theta)U(\omega,\theta)=U(\Theta\omega,\theta+a(g))
 \quad\text{for }m\text{-almost every }\theta
\]
for every $g\in\supp\nu_\omega$ and every $\omega$ in one full-measure
Borel set.  The exceptional phase set may depend on $(\omega,g)$.
  The cocycle
is environment-wise strongly irreducible if no such family exists.
\end{definition}
\begin{corollary}[Strong irreducibility implies stationary-value uniqueness]
\label[corollary]{corstrong-SVU}
Under \eqref{eq:log-integrability} and mixed-base ergodicity,
environment-wise strong irreducibility implies relative irreducibility,
and hence \SVU.
\end{corollary}
\begin{proof}
By the support-exactification lemma in Appendix~\ref{appsupport},
any relative invariant subbundle has support-level invariance with
exactly the preceding quantifiers.  Such a bundle is nonzero and
proper by definition, so its singleton family would contradict strong
irreducibility.  Apply \cref{thm:QI-SVU-equivalence}.
\end{proof}

\subsection{Uniform expected growth}
\label{subsecexpected}
The extremal conditional growth rates converge to the largest and
smallest stationary values. We first identify their integral limits
using \cref{lem:integrated-envelopes}, and then obtain almost-sure and
$L^1$ convergence from subadditivity and uniform integrability.

Set
\[
 a_n(\omega)=\max_{x\in Q}H_n(\omega,x),
 \qquad
 b_n(\omega)=\min_{x\in Q}H_n(\omega,x),
\]
and let
\[
 \alpha_-=
 \min_{\eta^\bullet\in\cS_Q(\nu^\bullet)}
          \alpha_{\nu^\bullet}(\eta^\bullet).
\]
For $\omega\in\Omega_*$, define
\[
 \Delta_n^{\mathrm{op}}(\omega)
 =\sup_{\theta\in\TT^d}
   \left|\frac1n\EE_\omega\log\|A_n(\xi,\theta)\|-L_1\right|.
\]

\begin{theorem}[Conditional expected-growth envelopes]
\label[theorem]{thmexpected}
Assume \eqref{eq:log-integrability} and mixed-base ergodicity. 
Then
\[
 \frac{a_n}{n}\longrightarrow L_1,
 \qquad
 \frac{b_n}{n}\longrightarrow\alpha_-
 \qquad P\text{-almost surely and in }L^1(P).
\]
In particular, \SVU\ holds if and only if
$\|\Delta_n\|_{L^1(P)}\to0$. Under \SVU, one also has
\begin{equation}\label{eqexpected-limit}
 \Delta_n\longrightarrow0
 \qquad P\text{-almost surely}.
\end{equation}

For every $n\geq1$ and $\omega\in\Omega_*$,
\[
 \Delta_n^{\mathrm{op}}(\omega)
 \leq\Delta_n(\omega)+\frac{C_r}{n},
\]
where the finite constant
\(C_r\)
depends only on $r$.
\end{theorem}

\begin{proof}
\emph{Integral limits.}
By \eqref{eqH-additive},
$H_n=S_n\psi_{\nu^\bullet}$. Apply
\cref{lem:integrated-envelopes} with
$\varphi=\psi_{\nu^\bullet}\in L^1(\Omega;C(Q))$.
The Furstenberg formula in \cref{thmstationary}(iii) identifies the
maximal stationary value with $L_1$, so
\begin{equation}\label{eqBN-envelopes}
 \lim_{n\to\infty}\frac1n\int_\Omega a_n\dd P=L_1,
 \qquad
 \lim_{n\to\infty}\frac1n\int_\Omega b_n\dd P=\alpha_-.
\end{equation}

\emph{Almost-sure and $L^1$ convergence.}
Conditional independence of the future samples gives
\[
 H_{n+k}(\omega,x)
 =H_n(\omega,x)
  +\EE_\omega H_k(\Theta^n\omega,X_n^x).
\]
Taking the maximum and minimum over $x$ yields
\[
 a_{n+k}\leq a_n+a_k\circ\Theta^n,
 \qquad
 b_{n+k}\geq b_n+b_k\circ\Theta^n.
\]
Also, \eqref{eqlog-envelope} gives
\[
 |a_n|,\ |b_n|
 \leq\sum_{j=0}^{n-1}h\circ\Theta^j,
 \qquad h\in L^1(P),\quad h\geq0.
\]
Kingman's subadditive ergodic theorem, applied to $a_n$ and $-b_n$,
therefore gives finite almost-sure limits for $a_n/n$ and $b_n/n$.
These limits are constant because $\Theta$ is ergodic.

To obtain $L^1$ convergence, put
\[
 V_n=\frac1n\sum_{j=0}^{n-1}h\circ\Theta^j,
 \qquad
 T_{n,R}=\frac1n\sum_{j=0}^{n-1}(h-R)_+\circ\Theta^j.
\]
For $R,T>0$, one has $V_n\leq R+T_{n,R}$ and, by invariance of $P$,
\[
 \int V_n\dd P=\|h\|_1,
 \qquad
 \int T_{n,R}\dd P=\int(h-R)_+\dd P.
\]
Consequently, Markov's inequality gives
\[
\begin{aligned}
 \sup_{n\geq1}\int V_n\one_{\{V_n>T\}}\dd P
 &\leq\sup_{n\geq1}
       \left(RP\{V_n>T\}+\int T_{n,R}\dd P\right)\\
 &\leq\frac{R\|h\|_1}{T}+\int(h-R)_+\dd P.
\end{aligned}
\]
First choose $R$ large and then $T$ large. This proves uniform
integrability of $(V_n)$ and hence of $(a_n/n)$ and $(b_n/n)$.
Their almost-sure convergence is therefore also convergence in
$L^1(P)$. The integral limits in \eqref{eqBN-envelopes} identify the
two limits as $L_1$ and $\alpha_-$, respectively.

\emph{Stationary-value uniqueness.}
Under \SVU, $\alpha_-=L_1$. Since
\[
 \Delta_n
 =\max\left\{
       \left|\frac{a_n}{n}-L_1\right|,
       \left|\frac{b_n}{n}-L_1\right|
       \right\},
\]
the two envelope limits imply $\|\Delta_n\|_{L^1(P)}\to0$ and \eqref{eqexpected-limit}.
Conversely, for every stationary family, \eqref{eqH-additive} and
the environment marginal of $\lambda^\eta$ give
\[
 \left|\alpha_{\nu^\bullet}(\eta^\bullet)-L_1\right|
 =\left|\int_{\Omega\times Q}
        \left(\frac{H_n}{n}-L_1\right)\dd\lambda^\eta\right|
 \leq\int_\Omega\Delta_n\dd P.
\]
Thus $\|\Delta_n\|_{L^1(P)}\to0$ implies \SVU.

\emph{Operator norms.}
Put
$G_n(\omega,\theta)=\EE_\omega\log\|A_n(\xi,\theta)\|$.
The calculation in the proof of \cref{thmstationary} gives
\eqref{eqprojective-comparison} with the stated choice of $C_r$.
For every $\omega\in\Omega_*$, integration of that comparison under
$\PP_\omega$ is justified by \eqref{eqlog-envelope}. Together with
$\|A_np\|\leq\|A_n\|\|p\|$, it yields
\[
 H_n(\omega,\theta,\hat p)
 \leq G_n(\omega,\theta)
 \leq\int_{\PP(\RR^r)}
           H_n(\omega,\theta,\hat q)\dd m_{\PP}(\hat q)+C_r.
\]
Hence $b_n(\omega)\leq G_n(\omega,\theta)\leq a_n(\omega)+C_r$,
so
\[
 L_1-\Delta_n(\omega)
 \leq\frac{G_n(\omega,\theta)}{n}
 \leq L_1+\Delta_n(\omega)+\frac{C_r}{n}.
\]
Taking the supremum over $\theta$ proves the last assertion.
\end{proof}

\subsection{Qualitative continuity of the top Lyapunov exponent}
\label{subsecqualcontinuity}
Under integrated Wasserstein convergence and uniform logarithmic tail
control, $L_1$ is continuous at a reference SVU kernel, provided all
mixed bases are ergodic (\cref{thmqualcontinuity}). The proof uses
fixed-time continuity of conditional expected growth and imposes no
SVU assumption on the approximating kernels.

Choose a bounded translation-invariant compatible torus metric and
use the bounded product metric
\begin{equation}\label{eqambient-metric}
 d_\cG(g,\widetilde g)
 =d_{\TT^d}(a(g),a(\widetilde g))
            +\min\{1,\|A(g)-A(\widetilde g)\|_\infty\}.
\end{equation}
Let $W_1$ be its Wasserstein distance.  For kernels modulo $P$-null sets, define
\begin{equation}\label{eqW-integrated}
 \W(\nu^\bullet,\widetilde\nu^\bullet)
 =\int_\Omega W_1(\nu_\omega,\widetilde\nu_\omega)\dd P(\omega).
\end{equation}
The bounded compatible cost makes $W_1$ Borel and metrizes weak
convergence. Thus $\W$ convergence is conditional weak convergence
in $P$-probability, a stronger requirement than weak convergence of
the joint laws. The tail hypothesis below supplies uniform
integrability for truncation; bounded first logarithmic moments alone
do not suffice. 

\begin{lemma}[Fixed-time stability with logarithmic tail control]
\label[lemma]{lemqualfinite}
Let $\nu^{(j),\bullet}$, $j\ge0$, satisfy
\eqref{eq:log-integrability}, with $j=0$ denoting $\nu^\bullet$.
Assume $\W(\nu^{(j),\bullet},\nu^\bullet)\to0$ and
\begin{equation}\label{equniform-log-tails}
 \lim_{R\to\infty}\sup_{j\ge0}
 \int_\Omega\int_\cG\Lambda_A(g)\one_{\{\Lambda_A(g)>R\}}
             \dd\nu^{(j)}_\omega(g)\dd P(\omega)=0.
\end{equation}
For each fixed $n$,
\[
 \|H_n^{\nu^{(j),\bullet}}-H_n^{\nu^\bullet}\|_{L^1(P;C(Q))}
 \longrightarrow0.
\]
\end{lemma}
\begin{proof}
From any subsequence choose a further subsequence for which the
conditional laws converge weakly for almost every environment.
Intersecting its finitely many shifts gives, for fixed $n$, weak
convergence of the product laws on $\cG^n$ at almost every $\omega$.
The finite-time integrand $Z_n(\mathbf g;x)$ is continuous in
$(\mathbf g,x)$. Clip it to $[-R,R]$. At a fixed good environment,
the convergent product probabilities form a tight family. On a
compact set carrying all but an arbitrarily small amount of their
mass, the family of clipped integrands indexed by compact $Q$ is
compact in the uniform norm. A finite net and weak convergence
therefore give convergence of the integrals uniformly in $x$.
The bound $2R$ permits integration over $P$.

The clipping error is dominated uniformly in $x$ by
$V\one_{\{V>R\}}$, where
$V=(r-1)\sum_{i<n}\Lambda_A(g_i)$. These variables, with varying
$j$ and joint laws $P(\dd\omega)\bigotimes_{i<n}
\nu^{(j)}_{\Theta^i\omega}(\dd g_i)$, are uniformly integrable.
Indeed, for nonnegative $a_i$,
\[
 \left(\sum_{i<n}a_i\right)\one_{\{\sum a_i>T\}}
 \le2\sum_{i<n}a_i\one_{\{a_i>T/(2n)\}}.
\]
Apply this with $a_i=(r-1)\Lambda_A(g_i)$, integrate, and use
$\Theta_*P=P$ and \eqref{equniform-log-tails}. Letting $R$ increase
proves the lemma along the further subsequence, and the subsequence
criterion gives full convergence.
\end{proof}

\begin{theorem}[Finite-scale stability and qualitative continuity]
\label[theorem]{thmqualcontinuity}
For any two logarithmically integrable kernels with ergodic mixed
bases, and every $n\ge1$,
\begin{equation}\label{eqfinite-qualitative}
\begin{split}
 |L_1(\widetilde\nu^\bullet)-L_1(\nu^\bullet)|
 \le \|\Delta_n^{\nu^\bullet}\|_1
   +\frac1n\|H_n^{\widetilde\nu^\bullet}-H_n^{\nu^\bullet}\|_{L^1(P;C(Q))}.
\end{split}
\end{equation}
Consequently, under the convergence and tail assumptions of
\cref{lemqualfinite}, if all mixed bases are ergodic and only the
limit $\nu^\bullet$ is assumed to have \SVU, then
\begin{equation}\label{eqqualcontinuity}
 L_1(\nu^{(j),\bullet})\longrightarrow L_1(\nu^\bullet).
\end{equation}
\end{theorem}
\begin{proof}
Take a maximizing stationary family $\widetilde\eta^\bullet$ for
$\widetilde\nu^\bullet$. The stationary identity gives
$nL_1(\widetilde\nu^\bullet)
=\int H_n^{\widetilde\nu^\bullet}\dd\lambda^{\widetilde\eta}$.
Add and subtract $H_n^{\nu^\bullet}$ and use that
$\lambda^{\widetilde\eta}$ has environment marginal $P$ to obtain
\eqref{eqfinite-qualitative}. For continuity, first choose $n$ so that
$\|\Delta_n^{\nu^\bullet}\|_1$ is small by \cref{thmexpected}, and
then use \cref{lemqualfinite} at this fixed $n$.
Under a common matrix and inverse bound, $\Lambda_A$ is bounded, so
the tail condition is automatic.
\end{proof}

\section{Quenched concentration in ergodic environments}
\label{secquenched}
SVU yields eventual quenched exponential concentration in an arbitrary
ergodic environment, with deterministic rates and an environment-dependent
starting scale (\cref{thmquenched}). The renewal example shows why this
does not imply annealed exponential decay. The forward estimates assume
that, for some $M\ge1$,
\begin{equation}\label{eqbounded-growth}
 \sup_{x\in Q}|\Psi(g,x)|\le M
 \quad\text{for }P(\dd\omega)\nu_\omega(\dd g)
                         \text{-almost every }(\omega,g).
\end{equation}
It implies $|H_n|\le Mn$ and $|L_1|\le M$.
Intersecting the relevant environment set with all its integer shifts
gives an invariant Borel set $\Omega_b$ of full measure where these
bounds hold along $\PP_\omega$-almost every sampled path. 
For $\varepsilon>0$, define
\begin{equation}\label{eqq-n}
 q_n(\omega,\varepsilon)
 =\sup_{x\in Q}\PP_\omega\{|Z_n^\omega(x)/n-L_1|>\varepsilon\}.
\end{equation}
Pathwise continuity and Fatou's lemma make the strict-event probability
lower semicontinuous in $x$. Its supremum over a countable dense subset
therefore equals $q_n$, proving measurability.

\subsection{Block martingale concentration}
The sliding-block estimate separates directional growth into an
environmental average of block drift and an exponentially controlled
martingale term. Its one-sided form will also apply to exterior growth.

\begin{lemma}[Sliding-block concentration]\label[lemma]{lemblock}
Assume \eqref{eqbounded-growth}. For $\lambda\in[-M,M]$ and
$\ell\ge1$, put
\[
 \Delta_\ell^\lambda(\omega)
 =\sup_x|H_\ell(\omega,x)/\ell-\lambda|,
 \qquad
 u_\ell^\lambda(\omega)
 =\bigl(\sup_xH_\ell(\omega,x)/\ell-\lambda\bigr)_+.
\]
For every $\omega\in\Omega_b$, $n\ge\ell$, and $t>0$,
\begin{equation}\label{eqblock-two-sided}
\begin{split}
 \sup_x\PP_\omega\left\{|Z_n^\omega(x)/n-\lambda|>
 \frac1n\sum_{k=0}^{n-1}\Delta_\ell^\lambda(\Theta^k\omega)
       +\frac{2M(\ell-1)}n+t\right\}
 \le 2\ell\exp\left(-\frac{nt^2}{8M^2\ell}\right).
\end{split}
\end{equation}
The corresponding upper-tail estimate is
\begin{equation}\label{eqblock-one-sided}
\begin{split}
 \sup_x\PP_\omega\left\{Z_n^\omega(x)/n-\lambda>
 \frac1n\sum_{k=0}^{n-1}u_\ell^\lambda(\Theta^k\omega)
       +\frac{2M(\ell-1)}n+t\right\}
 \le \ell\exp\left(-\frac{nt^2}{8M^2\ell}\right).
\end{split}
\end{equation}
\end{lemma}
\begin{proof}
Fix $\omega,x$ and complete
$\mathcal F_j^\omega=\sigma(\xi_0,\ldots,\xi_{j-1})$ under
$\PP_\omega$. Write $z_j=\Psi(\xi_j,X_j^x)$ and, for
$0\le k\le n-\ell$, put
\[
 S_{k,\ell}=\sum_{j=k}^{k+\ell-1}z_j,
 \qquad D_k=S_{k,\ell}-H_\ell(\Theta^k\omega,X_k^x).
\]
Conditional independence of the future samples gives
$\EE_\omega[S_{k,\ell}\mid\mathcal F_k^\omega]
=H_\ell(\Theta^k\omega,X_k^x)$. Thus $D_k$ is
$\mathcal F_{k+\ell}^\omega$-measurable, centered conditional on
$\mathcal F_k^\omega$, and $|D_k|\le2M\ell$.

For $a=0,\ldots,\ell-1$, let
\[
 N_a=\max\{0,\lfloor(n-a)/\ell\rfloor\},
 \qquad T_a=\sum_{q=0}^{N_a-1}D_{a+q\ell}.
\]
These indices partition $\{0,\ldots,n-\ell\}$, so
$\sum_aN_a=n-\ell+1$ and $N_a\le n/\ell$.
Along each residue class, the summands are martingale differences
for the filtration $\mathcal F_{a+q\ell}^\omega$.
The conditional Hoeffding bound \cite[inequality~(4.16)]{Hoeffding1963}
$\EE[e^{sD_{a+q\ell}}\mid\mathcal F_{a+q\ell}^\omega]
\le\exp(s^2(2M\ell)^2/2)$ gives, for $N_a>0$,
\[
 \PP_\omega\{|T_a|>v\}
 \le2\exp\left(-\frac{v^2}{8M^2\ell^2N_a}\right).
\]
For $N_a=0$ the event is empty. Therefore a union bound gives
\begin{equation}\label{eqall-residues}
 \PP_\omega\left\{\left|\sum_{k=0}^{n-\ell}D_k\right|>n\ell t\right\}
 \le2\ell\exp\left(-\frac{nt^2}{8M^2\ell}\right).
\end{equation}
No independence between the different residue classes is used.
For an upper-tail estimate the same argument has prefactor $\ell$.

Let $w_j$ count the length-$\ell$ blocks lying inside
$\{0,\ldots,n-1\}$ that contain $j$. Then $0\le w_j\le\ell$ and
\[
 \sum_{j<n}(\ell-w_j)
 =n\ell-\ell(n-\ell+1)=\ell(\ell-1).
\]
Consequently
\begin{equation}\label{eqsliding-identity}
 \ell(Z_n-n\lambda)
 =\sum_{k=0}^{n-\ell}D_k
  +\sum_{k=0}^{n-\ell}
     \bigl(H_\ell(\Theta^k\omega,X_k^x)-\ell\lambda\bigr)
  +\sum_{j<n}(\ell-w_j)(z_j-\lambda).
\end{equation}
The last sum has absolute value at most
$2M\ell(\ell-1)$. Divide by $n\ell$, bound each middle summand by
$\ell\Delta_\ell^\lambda(\Theta^k\omega)$, and extend the nonnegative
environmental sum to $n-1$. This and \eqref{eqall-residues} prove
\eqref{eqblock-two-sided}. Keeping signs and bounding the middle
summands above by $\ell u_\ell^\lambda$ proves
\eqref{eqblock-one-sided}. All constants are independent of $x$.
\end{proof}

\Needspace{12\baselineskip}
\subsection{Quenched exponential concentration}
Under SVU, exponential quenched rates hold after an almost surely
finite environmental cutoff. Only Birkhoff convergence, not a mixing
rate, is needed for that cutoff.

\begin{theorem}[Eventual quenched exponential concentration]
\label[theorem]{thmquenched}
Assume mixed-base ergodicity, \SVU, and \eqref{eqbounded-growth}.
For every $\varepsilon>0$ there are deterministic
$C_\varepsilon,c_\varepsilon>0$ and a measurable almost surely finite
$N_\varepsilon(\omega)$ such that
\begin{equation}\label{eqquenched-exponential}
 q_n(\omega,\varepsilon)\le C_\varepsilon e^{-c_\varepsilon n}
 \quad(n\ge N_\varepsilon(\omega))
\end{equation}
for $P$-almost every $\omega$.
\end{theorem}
\begin{proof}
By \cref{thmexpected}, choose $\ell$ with
$\int\Delta_\ell\dd P\le\varepsilon/8$. Birkhoff's theorem gives
$N_1(\omega)<\infty$ almost surely such that
$n^{-1}\sum_{j<n}\Delta_\ell(\Theta^j\omega)\le\varepsilon/4$
for all $n\ge N_1(\omega)$. One may choose the least such eventual
cutoff, which is measurable since all these inequalities are
measurable and indexed by integers. Take
\[
 N_\varepsilon(\omega)
 =\max\{N_1(\omega),\ell,\lceil8M\ell/\varepsilon\rceil\}.
\]
Apply \eqref{eqblock-two-sided} with $\lambda=L_1$ and
$t=\varepsilon/2$. The full threshold is at most $\varepsilon$.
Thus $C_\varepsilon=2\ell$ and
$c_\varepsilon=\varepsilon^2/(32M^2\ell)$ work.
\end{proof}

Conversely, integrated convergence of the conditional deviation
probabilities forces SVU, even without bounded growth.

\begin{proposition}[A logarithmically integrable converse]
\label[proposition]{propconverse}
Assume mixed-base ergodicity and \eqref{eq:log-integrability}, but
not necessarily bounded growth. If, for every $\varepsilon>0$,
\[
 \int q_n(\omega,\varepsilon)\dd P(\omega)\longrightarrow0,
\]
then the kernel has \SVU. In particular this holds if
$q_n(\omega,\varepsilon)\to0$ almost surely for every
$\varepsilon>0$.
\end{proposition}
\begin{proof}
Fix a stationary family $\eta^\bullet$. Under the joint law
$\rho_\eta=P(\dd\omega)\eta_\omega(\dd x)\PP_\omega(\dd\xi)$,
the assumed bound gives $Z_n/n\to L_1$ in probability. The variables
are dominated by
$V_n=(r-1)n^{-1}\sum_{j<n}\Lambda_A(\xi_j)$.
For $U=(r-1)\Lambda_A(\xi_0)$, all the shifted summands have the
same law under $\rho_\eta$. Splitting them into a part bounded by
$R$ and a remainder gives the uniform-integrability estimate
\[
 \sup_n\int V_n\one_{\{V_n>T\}}\dd\rho_\eta
 \le \frac{R\,\EE_{\rho_\eta}U}{T}
                 +\EE_{\rho_\eta}(U-R)_+ \quad(T>2R).
\]
Therefore convergence in probability is convergence in $L^1$.
The expectation of $Z_n/n$ equals
$\alpha_{\nu^\bullet}(\eta^\bullet)$ for every $n$, by
\eqref{eqH-additive}. It must equal $L_1$.
Finally $0\le q_n\le1$ allows dominated convergence for the last
assertion.
\end{proof}

For later use define
\begin{equation}\label{eqq-op}
 q_n^{\mathrm{op}}(\omega,\varepsilon)
 =\sup_\theta\PP_\omega\{|n^{-1}\log\|A_n(\xi,\theta)\|-L_1|
                                             >\varepsilon\}.
\end{equation}
If $\log r/(2n)\le\varepsilon/2$, then
\begin{equation}\label{eqop-comparison}
 q_n^{\mathrm{op}}(\omega,\varepsilon)
 \le(r+1)q_n(\omega,\varepsilon/2).
\end{equation}
Indeed, for an orthonormal basis $e_i$,
$\|Te_1\|\le\|T\|\le\sqrt r\max_i\|Te_i\|$.
The lower norm tail is contained in the tail for $e_1$, and the upper
norm tail in the union of the $r$ upper directional tails.
This transfers both quenched and integrated estimates to operator
norms without a separate structural assumption.

\subsection{Environmental decomposition and obstruction}
Control of the environmental Birkhoff term below transfers the
quenched estimate to an annealed one. The renewal example shows that
general ergodicity alone cannot supply annealed exponential rates.

\begin{proposition}[Environmental decomposition]
\label[proposition]{propenvironment-decomposition}
Under mixed-base ergodicity and \eqref{eqbounded-growth}, for every
$\ell\ge1$, $\varepsilon>0$, and
$n\ge\max\{\ell,\lceil8M\ell/\varepsilon\rceil\}$,
\begin{equation}\label{eqenvironment-decomposition}
\begin{split}
 \int q_n(\omega,\varepsilon)\dd P(\omega)
 &\le2\ell\exp\left(-\frac{\varepsilon^2n}{32M^2\ell}\right)\\
 &\quad+P\left\{\frac1n\sum_{j<n}\Delta_\ell(\Theta^j\omega)
                                      >\frac\varepsilon4\right\}.
\end{split}
\end{equation}
\end{proposition}
\begin{proof}
Outside the environmental event, use \cref{lemblock} with
$\lambda=L_1$ and $t=\varepsilon/2$. On that event use $q_n\le1$,
and integrate. 
\end{proof}

\begin{proposition}[A renewal obstruction]\label[proposition]{proprenewal-obstruction}
There is an invertible ergodic environment and a compactly supported
two-dimensional moving-law mixed cocycle with \SVU\ and $L_1>0$
for which, at $\varepsilon=L_1/2$,
\[
 \int q_n(\omega,\varepsilon)\dd P(\omega)\ge c n^{-2}
 \quad(n\ge1),\qquad
 \esssup_\omega q_n(\omega,\varepsilon)=1.
\]
Thus neither annealed exponential decay nor a deterministic uniform
quenched starting scale follows from general ergodicity and \SVU.
\end{proposition}
\begin{proof}
Let a Markov chain $J$ on $\{0,1,\ldots\}$ count down from $j\ge1$
to $j-1$. From $0$, it jumps to an independent $K_*$ with
$P(K_*=k)=c_*k^{-4}$, $k\ge1$. Its mean return time to zero is
$m_*=1+\EE K_*<\infty$. The stationary probabilities are
\[
 \pi_0=m_*^{-1},\qquad
 \pi_j=m_*^{-1}P(K_*\ge j)\quad(j\ge1).
\]
They sum to one and satisfy the balance equations. The chain is
irreducible and positive recurrent. Return lengths $2$ and $3$
make it aperiodic. Its stationary two-sided shift
$(\Omega,P,\Theta)$ is therefore invertible and ergodic.
Fix $a_*>0$ and put $a(\omega)=a_*\one_{\{J_0=0\}}$.
Conditionally on $\omega$, sample independently a Haar translation
$\alpha$ and a Haar rotation $R\in\mathrm{SO}(2)$, and use
\[
 g=(\alpha,D(a(\omega))R),\qquad
 D(a)=\operatorname{diag}(e^a,e^{-a}),
\]
with a constant matrix field in the phase. These laws share a compact
support and satisfy \eqref{eqbounded-growth} with
$M=\max\{1,a_*\}$. The mixed base is ergodic by
\cref{lemtranslation-ergodicity}.

For every starting $x$, the projective transition sends $x$ to
$\zeta_{a(\omega)}=m\otimes D(a(\omega))_*m_{\PP}$.
It follows that the unique stationary family is
$\eta_\omega=\zeta_{a(\Theta^{-1}\omega)}$.
Moreover,
\begin{equation}\label{eqrenewal-one-step}
 \int\log\|D(a)u\|\dd m_{\PP}(\hat u)=\log\cosh a.
\end{equation}
For $u=(\cos t,\sin t)$, the logarithmic norm is
$\log\cosh a+\log|1+\tanh(a)e^{-2it}|$. The latter term has
zero circle average, by its absolutely convergent power series when
$a$ is finite (and trivially when $a=0$). Thus \SVU\ holds and
$L_1=\pi_0\log\cosh a_*>0$ by the Furstenberg formula.

The conditional drift of each $z_j$ is
$b(\Theta^j\omega)=\log\cosh a(\Theta^j\omega)$, independently of
$X_j$. Since $z_j$ has conditional range of length at most $2a_*$,
conditional Hoeffding and iteration give the explicit bound
\[
 \sup_x\PP_\omega\left\{
 \left|Z_n^\omega(x)-\sum_{j<n}b(\Theta^j\omega)\right|>nt\right\}
 \le2\exp\left(-\frac{nt^2}{2a_*^2}\right).
\]
On $E_n=\{J_0\ge n\}$, all of the first $n$ matrices are rotations,
so $Z_n^\omega(x)=0$ for every $x$ and every allowed sampled path.
Hence $q_n(\omega,L_1/2)=1$ on $E_n$. Finally,
\[
 P(E_n)=\frac{c_*}{m_*}
       \sum_{k\ge n}(k-n+1)k^{-4}\asymp n^{-2}.
\]
The last comparison follows by summing over $2n\le k\le3n$ for
the lower bound and using $k-n+1\le k$ for the upper bound.
It proves both assertions.
\end{proof}

\section{Uniform exponential deviations in finite-state environments}
\label{secfinite-state}

For a fixed finite-state irreducible aperiodic Markov environment,
annealed exponential deviation estimates are uniform over nearby
sampling kernels at each fixed tolerance, and a common quenched cutoff
has an exponential tail. It suffices to assume SVU at
the reference. Fixed-time continuity transfers the small reference block
error to nearby kernels, and word-chain concentration controls the
environmental term in \eqref{eqenvironment-decomposition}.

\subsection{Finite-state kernels and fixed-time continuity}
\label{subsecfinite-framework}
Fix a finite set $S$ and an irreducible aperiodic transition matrix
$\Pi=(\Pi_{st})_{s,t\in S}$, with stationary probability $\pi$.
Let $(\Omega,P,\Theta)$ be its stationary two-sided Markov shift on
the space of allowed sequences, equipped with the product topology.
Since $\Pi$ is irreducible and aperiodic, this shift is mixing and
hence ergodic. 
Let $\boldsymbol\nu=(\nu_s)_{s\in S}$ be a probability kernel
from $S$ to $\cG$. Set
\begin{equation}\label{eqstate-dependent-law}
\nu_\omega=\nu_{\omega_0},
\qquad
\PP_{\boldsymbol\nu,\omega}
=\bigotimes_{j\ge0}\nu_{\omega_j}.
\end{equation}
Throughout Sections~\ref{secfinite-state}--\ref{secschrodinger},
$(\Pi,\pi,P)$ is fixed and only the sampling kernel $\boldsymbol\nu$ varies.

Fix $B\ge e$ and set
\begin{equation}\label{eqcompact-bounds}
 \cG_B=
 \left\{g\in\cG:
 \sup_{\theta\in\TT^d}
 \max\bigl\{\|A(g)(\theta)\|,\|A(g)(\theta)^{-1}\|\bigr\}
 \le B\right\},
 \qquad M=\log B.
\end{equation}
The class $\cG_B$ is a closed Polish subspace of $\cG$, although it
is generally noncompact. In this section all kernels belong to
$\Prob(\cG_B)^S$.  The common
bound \cref{eqcompact-bounds} allows us to take $\Omega_b=\Omega$ in
\cref{lemblock}, simultaneously for all these kernels.

The ambient integrated metric from \eqref{eqW-integrated} becomes
\begin{equation}\label{eqkernel-distance}
 \W(\boldsymbol\nu,\widetilde{\boldsymbol\nu})
 =\sum_{s\in S}\pi_sW_1(\nu_s,\widetilde\nu_s),
\end{equation}
where $W_1$ uses the restriction of \eqref{eqambient-metric} to
$\cG_B$. The auxiliary quantity
$d_{\max}=\max_sW_1(\nu_s,\widetilde\nu_s)$ satisfies
\begin{equation}\label{eqmax-equivalence}
 \W\le d_{\max}\le\pi_{\min}^{-1}\W,
 \qquad \pi_{\min}=\min_s\pi_s>0.
\end{equation}
Set
\begin{equation}\label{eqergodic-kernel-class}
 \mathcal E_B=
 \left\{\boldsymbol\nu\in\Prob(\cG_B)^S:
 (F,\rho_{\boldsymbol\nu})\text{ is ergodic}\right\}.
\end{equation}
Neighborhoods are relative to $\mathcal E_B$, which need not be open
in $\Prob(\cG_B)^S$. Superscripts on $H_n^{\boldsymbol\nu}$,
$\Delta_n^{\boldsymbol\nu}$, and $q_n^{\boldsymbol\nu}$ specify the
kernel. The latter two quantities are centered at that kernel's
exponent $L_1(\boldsymbol\nu)$.

\begin{lemma}[Fixed-time kernel continuity]
\label[lemma]{lemfinite-continuity}
For every fixed $\ell\ge1$,
$H_\ell^{\boldsymbol\nu}(\omega,x)$ depends on $\omega$ only through
the finite block $(\omega_0,\ldots,\omega_{\ell-1})$. The map
\(
\boldsymbol\nu\longmapsto H_\ell^{\boldsymbol\nu}
\)
is continuous from $(\Prob(\cG_B)^S,\W)$ to
$C(\Omega\times Q)$ endowed with the uniform norm. The same
finite-block dependence and continuity hold for the conditional growth
function of the exterior-square cocycle on its projective state space.

If $\boldsymbol\nu_0\in\mathcal E_B$ has \SVU, then, for every fixed
$\ell\ge1$,
\begin{equation}\label{eqdelta-neighborhood}
 \|\Delta_\ell^{\boldsymbol\nu}
      -\Delta_\ell^{\boldsymbol\nu_0}\|_\infty
 \longrightarrow0
 \qquad\text{as }\boldsymbol\nu\to\boldsymbol\nu_0
 \text{ in }\W\text{ within }\mathcal E_B.
\end{equation}
\end{lemma}

\begin{proof}
For an allowed word $w=(s_0,\ldots,s_{\ell-1})$,
\[
 H_\ell^{\boldsymbol\nu}(w,x)
 =\int_{\cG_B^\ell}Z_\ell(\mathbf g;x)
                   \prod_{i=0}^{\ell-1}d\nu_{s_i}(g_i).
\]
The integrand is jointly continuous and bounded by $\ell M$.
If $\boldsymbol\nu^{(j)}\to\boldsymbol\nu$ in $\W$, then
\eqref{eqmax-equivalence} gives weak convergence of each coordinate
law and hence of the product laws for $w$. The tightness and finite-net
argument in \cref{lemqualfinite}, now with a bounded integrand,
yields convergence uniformly in $x$. Explicitly, on a compact set
carrying mass at least $1-\delta$ for all these products, a finite
$\eta$-net for $\{Z_\ell(\cdot;x):x\in Q\}$ bounds the uniform
error by the maximum of finitely many weak-convergence errors plus
$2\eta+4\ell M\delta$. There are finitely many allowed words, so
convergence is also uniform in $\omega$. The same argument applies to
the exterior integrand, bounded by $2\ell M$.

Finally,
\[
 \|\Delta_\ell^{\boldsymbol\nu}
       -\Delta_\ell^{\boldsymbol\nu_0}\|_\infty
 \le\ell^{-1}\|H_\ell^{\boldsymbol\nu}
                  -H_\ell^{\boldsymbol\nu_0}\|_\infty
       +|L_1(\boldsymbol\nu)-L_1(\boldsymbol\nu_0)|.
\]
The first term tends to zero by the preceding argument. The second
tends to zero by \cref{thmqualcontinuity}, since $\Lambda_A\le M$
provides uniform logarithmic tails and the kernels lie in $\mathcal E_B$.
\end{proof}

\subsection{Environmental word chains and concentration}
\label{subsecword-chain}

For $\ell\ge1$, let
\[
\mathcal W_\ell
=
\left\{(w_0,\ldots,w_{\ell-1})\in S^\ell:
       \prod_{i=0}^{\ell-2}\Pi_{w_i,w_{i+1}}>0\right\},
\]
with an empty product interpreted as $1$. The process
$Y_j=(\omega_j,\ldots,\omega_{j+\ell-1})$ is a stationary Markov
chain on $\mathcal W_\ell$, with transition matrix
\begin{equation}\label{eqword-transition}
R_\ell(w,w')
=
\mathbf1_{\{(w_1,\ldots,w_{\ell-1})
                 =(w'_0,\ldots,w'_{\ell-2})\}}
\Pi_{w_{\ell-1},w'_{\ell-1}}.
\end{equation}
The overlap condition is vacuous when $\ell=1$.
Its stationary law, obtained as the distribution of $Y_0$ under $P$, is
\begin{equation}\label{eqword-stationary}
\pi_\ell(w)
=
\pi_{w_0}\prod_{i=0}^{\ell-2}\Pi_{w_i,w_{i+1}}.
\end{equation}
For $h:\mathcal W_\ell\to\RR$, write
$\pi_\ell h=\sum_{w\in\mathcal W_\ell}\pi_\ell(w)h(w)$.
We use the total-variation convention
\[
\|\mu-\zeta\|_{\mathrm{TV}}
=\sup_E|\mu(E)-\zeta(E)|.
\]

\begin{lemma}[Finite-window concentration]
\label[lemma]{lemword-concentration}
For every fixed $\ell\ge1$, the number
\begin{equation}\label{eqword-mixing-time}
T_\ell
=
1+\sum_{q=0}^{\infty}
\sup_{w\in\mathcal W_\ell}
\|R_\ell^q(w,\cdot)-\pi_\ell\|_{\mathrm{TV}}
\end{equation}
is finite. If $H_0>0$, $h:\mathcal W_\ell\to[0,H_0]$,
$\delta>0$, and $n\ge1$ satisfies $n\delta\ge4H_0T_\ell$, then
\begin{equation}\label{eqword-concentration}
P\left\{
\frac1n\sum_{j=0}^{n-1}h(Y_j)-\pi_\ell h>\delta
\right\}
\le
\exp\left(-\frac{n\delta^2}{32H_0^2T_\ell^2}\right).
\end{equation}
The bound is uniform over $h$ with this range restriction, and
$T_\ell$ depends only on $\ell$ and the fixed environment chain.
\end{lemma}

\begin{proof}
For $q\ge\ell-1$, the Markov property gives
\[
R_\ell^q(w,w')
=
(\Pi^{q-\ell+1})_{w_{\ell-1},w'_0}
\prod_{i=0}^{\ell-2}\Pi_{w'_i,w'_{i+1}}.
\]
Primitivity of $\Pi$ gives $p_0\ge1$ such that $\Pi^{p_0}>0$.
Thus $R_\ell^{q_0}>0$ for $q_0=p_0+\ell-1$, and for some
$0<\beta<1$,
\[
 R_\ell^{q_0}(w,\cdot)\ge\beta\upsilon_\ell(\cdot),
\]
where $\upsilon_\ell$ is uniform on $\mathcal W_\ell$.
Writing $R_\ell^{q_0}=\beta\upsilon_\ell+(1-\beta)\mathsf Q$
and using total-variation contraction of Markov kernels gives
\[
 \|\mu R_\ell^{q_0}-\zeta R_\ell^{q_0}\|_{\mathrm{TV}}
 \le(1-\beta)\|\mu-\zeta\|_{\mathrm{TV}}.
\]
Iteration and stationarity imply
\[
 \sup_w\|R_\ell^q(w,\cdot)-\pi_\ell\|_{\mathrm{TV}}
 \le(1-\beta)^{\lfloor q/q_0\rfloor},
\]
so $T_\ell<\infty$.

Set $f=h-\pi_\ell h$, so that $\pi_\ell f=0$. For each $q\ge0$,
\[
R_\ell^qf(w)=R_\ell^qh(w)-\pi_\ell h,
\qquad
|R_\ell^qf(w)|
\le H_0\|R_\ell^q(w,\cdot)-\pi_\ell\|_{\mathrm{TV}}.
\]
Here the range of $h$ and the stated TV convention give the second
inequality. The Poisson series
\[
u=\sum_{q=0}^{\infty}R_\ell^qf
\]
converges uniformly and satisfies $\|u\|_\infty\le H_0T_\ell$.
It solves $(I-R_\ell)u=f$.

Let $\mathcal F_j^Y=\sigma(Y_0,\ldots,Y_j)$ and define
\[
m_{j+1}=u(Y_{j+1})-(R_\ell u)(Y_j).
\]
These are $\mathcal F_{j+1}^Y$-measurable martingale differences
with $|m_{j+1}|\le2H_0T_\ell$.
The Poisson equation implies
\[
f(Y_j)=u(Y_j)-u(Y_{j+1})+m_{j+1},
\]
and hence
\begin{equation}\label{eqword-poisson}
\sum_{j=0}^{n-1}f(Y_j)
=
u(Y_0)-u(Y_n)+\sum_{j=0}^{n-1}m_{j+1}.
\end{equation}
The endpoint term has absolute value at most
$2H_0T_\ell\le n\delta/2$. Therefore one-sided Azuma--Hoeffding
\cite[Section~2]{Azuma1967} gives
\[
\begin{aligned}
P\left\{\sum_{j=0}^{n-1}f(Y_j)>n\delta\right\}
&\le
P\left\{\sum_{j=0}^{n-1}m_{j+1}>\frac{n\delta}{2}\right\}\\
&\le
\exp\left(-\frac{(n\delta/2)^2}{2n(2H_0T_\ell)^2}\right)
=
\exp\left(-\frac{n\delta^2}{32H_0^2T_\ell^2}\right).
\end{aligned}
\]
This is \eqref{eqword-concentration}.
\end{proof}

\subsection{Annealed exponential deviations uniform over nearby kernels}
\label{subsecuniform-deviations}

At each fixed tolerance, annealed exponential deviations hold on one
relative neighborhood of a reference SVU kernel. The neighborhood and
rate constants may depend on that tolerance, but not on the nearby
kernel.

\begin{proposition}[Uniform environmental deviations]
\label[proposition]{propuniform-environment}
Let $\boldsymbol\nu_0\in\mathcal E_B$ have \SVU. For every
$\varepsilon>0$, there exist an integer $\ell\ge1$ and a relative
neighborhood $\mathcal U_\varepsilon\subset\mathcal E_B$ of
$\boldsymbol\nu_0$ such that
\begin{equation}\label{equniform-small-block}
\int_\Omega\Delta_\ell^{\boldsymbol\nu_0}\,dP
\le\frac{\varepsilon}{16},\qquad
\|\Delta_\ell^{\boldsymbol\nu}
       -\Delta_\ell^{\boldsymbol\nu_0}\|_\infty
\le\frac{\varepsilon}{16},\qquad
\int_\Omega\Delta_\ell^{\boldsymbol\nu}\,dP
\le\frac{\varepsilon}{8}
\end{equation}
for every $\boldsymbol\nu\in\mathcal U_\varepsilon$. For every
such kernel and every integer
$n\ge\lceil64MT_\ell/\varepsilon\rceil$,
\begin{equation}\label{equniform-environment-ldt}
P\left\{
\frac1n\sum_{j=0}^{n-1}
\Delta_\ell^{\boldsymbol\nu}(\Theta^j\omega)
>
\frac{\varepsilon}{4}
\right\}
\le
\exp\left(-\frac{\varepsilon^2n}{8192M^2T_\ell^2}\right).
\end{equation}
\end{proposition}

\begin{proof}
Choose $\ell$ with
$\int\Delta_\ell^{\boldsymbol\nu_0}\,dP\le\varepsilon/16$
by \cref{thmexpected}, and then choose $\mathcal U_\varepsilon$
so that the uniform difference is at most $\varepsilon/16$ by
\cref{lemfinite-continuity}. Integration gives the third bound in
\eqref{equniform-small-block}.

For each nearby kernel, $\Delta_\ell^{\boldsymbol\nu}$ is a word
observable $h_{\boldsymbol\nu}$ with
$0\le h_{\boldsymbol\nu}\le2M$ and
$\pi_\ell h_{\boldsymbol\nu}\le\varepsilon/8$. Consequently,
\[
 \left\{\frac1n\sum_{j<n}
     \Delta_\ell^{\boldsymbol\nu}(\Theta^j\omega)>
     \frac\varepsilon4\right\}
 \subset
 \left\{\frac1n\sum_{j<n}h_{\boldsymbol\nu}(Y_j)
          -\pi_\ell h_{\boldsymbol\nu}>\frac\varepsilon8\right\}.
\]
Apply \cref{lemword-concentration} with $H_0=2M$ and
$\delta=\varepsilon/8$. Its scale restriction and exponent are exactly
those in \eqref{equniform-environment-ldt}.
\end{proof}

Combining \eqref{equniform-environment-ldt} with the sampling bound in
\eqref{eqenvironment-decomposition} gives the annealed estimate.

\begin{theorem}[Annealed exponential deviations uniform over nearby kernels]
\label[theorem]{thmuniform-ldt}
Fix the finite-state environment, $B\ge e$, and
$\boldsymbol\nu_0\in\mathcal E_B$ with \SVU. For every
$\varepsilon>0$, there exist a relative neighborhood
$\mathcal U_\varepsilon\subset\mathcal E_B$ of $\boldsymbol\nu_0$,
deterministic constants $C_\varepsilon,c_\varepsilon>0$, and an integer
$N_\varepsilon\ge1$ such that
\begin{equation}\label{equniform-annealed-ldt}
\sup_{\boldsymbol\nu\in\mathcal U_\varepsilon}
\int_\Omega q_n^{\boldsymbol\nu}(\omega,\varepsilon)\,dP(\omega)
\le
C_\varepsilon e^{-c_\varepsilon n}
\qquad(n\ge N_\varepsilon).
\end{equation}
The same conclusion holds
for $q_n^{\mathrm{op},\boldsymbol\nu}$, with possibly different
constants, starting scale, and neighborhood.
\end{theorem}

\begin{proof}
Choose $\ell$ and $\mathcal U_\varepsilon$ as in
\cref{propuniform-environment}. Combining
\cref{eqenvironment-decomposition,equniform-environment-ldt}, for every
$\boldsymbol\nu\in\mathcal U_\varepsilon$ we obtain
\[
\int_\Omega q_n^{\boldsymbol\nu}(\omega,\varepsilon)\,dP(\omega)
\le{}
2\ell\exp\left(-\frac{\varepsilon^2n}{32M^2\ell}\right)
+
\exp\left(-\frac{\varepsilon^2n}{8192M^2T_\ell^2}\right)
\]
whenever $n$ meets both scale restrictions. Thus we may take
\begin{equation}\label{equniform-constants}
\begin{aligned}
C_\varepsilon&=2\ell+1,\\
c_\varepsilon&=
\min\left\{
\frac{\varepsilon^2}{32M^2\ell},
\frac{\varepsilon^2}{8192M^2T_\ell^2}
\right\},\\
N_\varepsilon&=
\max\left\{
\ell,
\left\lceil\frac{8M\ell}{\varepsilon}\right\rceil,
\left\lceil\frac{64MT_\ell}{\varepsilon}\right\rceil
\right\}.
\end{aligned}
\end{equation}
These choices are independent of $\boldsymbol\nu$.

For operator norms, apply \eqref{eqop-comparison} to the directional
estimate at tolerance $\varepsilon/2$. Use
$\mathcal U_{\varepsilon/2}$, multiply the prefactor by $r+1$, and
take $n\ge\max\{N_{\varepsilon/2},\lceil(\log r)/\varepsilon\rceil\}$.
\end{proof}

\subsection{Common random cutoffs}
\label{subseccommon-cutoffs}

A single measurable quenched cutoff, with an exponential tail, works
for every kernel in the neighborhood. This requires the pointwise
reference-error comparison in \eqref{equniform-small-block}; the
uniform annealed bound alone does not give a common exceptional set.

\begin{corollary}[A common random starting scale]
\label[corollary]{corcommon-cutoff}
Under the hypotheses of \cref{thmuniform-ldt}, fix $\varepsilon>0$.
The neighborhood $\mathcal U_\varepsilon$ can be chosen so that there
exist deterministic constants $C,c,C',c'>0$, a measurable
$N:\Omega\to\mathbb N\cup\{\infty\}$, and a set
$\Omega_\varepsilon\subset\Omega$ of full $P$-measure with
$N(\omega)<\infty$ for every $\omega\in\Omega_\varepsilon$, such that
\begin{equation}\label{eqcommon-cutoff}
\begin{aligned}
q_n^{\boldsymbol\nu}(\omega,\varepsilon)
&\le Ce^{-cn}
&&\bigl(\boldsymbol\nu\in\mathcal U_\varepsilon,
        \ n\ge N(\omega)\bigr),\\
P\{N>n\}&\le C'e^{-c'n}
&&(n\ge1).
\end{aligned}
\end{equation}
Both this environmental set and the cutoff $N(\omega)$ are independent
of the nearby kernel, at the fixed tolerance $\varepsilon$.
\end{corollary}

\begin{proof}
Choose $\ell$ and $\mathcal U_\varepsilon$ as in
\cref{propuniform-environment}, and define
\[
E_n=
\left\{
\omega\in\Omega:
\frac1n\sum_{j=0}^{n-1}
\Delta_\ell^{\boldsymbol\nu_0}(\Theta^j\omega)
>
\frac{3\varepsilon}{16}
\right\}.
\]
The reference mean is at most $\varepsilon/16$, so
\cref{lemword-concentration} with $H_0=2M$ and
$\delta=\varepsilon/8$ gives
\[
P(E_n)\le e^{-c'n},
\qquad
c'=\frac{\varepsilon^2}{8192M^2T_\ell^2},
\qquad
n\ge\left\lceil\frac{64MT_\ell}{\varepsilon}\right\rceil.
\]
For every $\omega\notin E_n$, the pointwise kernel comparison in
\cref{equniform-small-block} gives simultaneously for all
$\boldsymbol\nu\in\mathcal U_\varepsilon$,
\[
\begin{aligned}
\frac1n\sum_{j=0}^{n-1}
\Delta_\ell^{\boldsymbol\nu}(\Theta^j\omega)
&\le
\frac1n\sum_{j=0}^{n-1}
\Delta_\ell^{\boldsymbol\nu_0}(\Theta^j\omega)
+\frac{\varepsilon}{16}\\
&\le\frac{3\varepsilon}{16}+\frac{\varepsilon}{16}
=\frac{\varepsilon}{4}.
\end{aligned}
\]

Choose a deterministic $N_0$ at least as large as the three quantities
defining $N_\varepsilon$ in \cref{equniform-constants}, and set
\[
N(\omega)=
\max\left\{
N_0,
1+\sup\{k\ge N_0:\omega\in E_k\}
\right\},
\]
where the inner supremum is $N_0-1$ for the empty set and $\infty$
for an unbounded set. This cutoff is measurable, being defined by
countably many measurable events. Since
$\sum_{k\ge N_0}P(E_k)<\infty$, the first Borel--Cantelli lemma gives
$P(\Omega_\varepsilon)=1$ for
$\Omega_\varepsilon=\{\omega:N(\omega)<\infty\}$.
For each integer $n\ge N_0$,
\[
\{N>n\}\subset\bigcup_{k\ge n}E_k,
\qquad
P\{N>n\}
\le\sum_{k\ge n}e^{-c'k}
=\frac{e^{-c'n}}{1-e^{-c'}}.
\]
Increasing $C'$ covers the finitely many integers $1\le n<N_0$.

For $\omega\in\Omega_\varepsilon$ and $n\ge N(\omega)$, the
environmental average is at most $\varepsilon/4$ simultaneously for
all nearby kernels. The definition of $N_0$ also makes the boundary
term at most $\varepsilon/4$. Thus \eqref{eqblock-two-sided} with
$t=\varepsilon/2$ gives \eqref{eqcommon-cutoff} with $C=2\ell$ and
$c=\varepsilon^2/(32M^2\ell)$ on the common full-measure set.
\end{proof}

\section{Quantitative stability and H\"older continuity}\label{secholder}

With the finite-state environment and current-state dependence
\eqref{eqstate-dependent-law} fixed, $L_1$ is locally two-point
H\"older continuous near a reference SVU kernel with a simple top
exponent (\cref{thmholder}). Top-norm deviations at scales $n$ and $2n$,
together with one-sided exterior-square upper deviations, provide
finite-scale singular-value separation and non-cancellation, allowing
the Avalanche Principle to approximate $L_1$ uniformly and exponentially
accurately by the corrected averages $\mathcal F_n$. The finite-time
kernel comparisons of \cref{propfinite-comparison} then yield the
H\"older bounds in $\Wsyn$ for fixed translation marginals and bounded
continuous matrix fields, and in $\W$ for varying translations under a
common phase H\"older bound, within the mixed-base-ergodic class.

\subsection{Finite-time perturbation estimates}
\label{secfinite-stability}
Retain $B\ge e$, $\cG_B$, and $M=\log B$ from
\eqref{eqcompact-bounds}. For
$\boldsymbol\nu\in\Prob(\cG_B)^S$, set
\[
 \mathcal L_n(\boldsymbol\nu)
 =\frac1n\int_\Omega\int_{\TT^d}
       \EE_{\boldsymbol\nu,\omega}\log\|A_n(\xi,\theta)\|
                                  \dd m(\theta)\dd P(\omega).
\]
For fixed probability marginals $\boldsymbol\mu=(\mu_s)_{s\in S}$,
write
\begin{equation}\label{eqfixed-marginal-class}
 \mathcal C_{\boldsymbol\mu}(B)
 =\{\boldsymbol\nu\in\Prob(\cG_B)^S:
                      a_*\nu_s=\mu_s\ \text{for all }s\}.
\end{equation}
For two laws with the same translation marginal, let
$\operatorname{Cpl}_{\mathrm{syn}}(\nu,\widetilde\nu)$ be the couplings
supported on $\{a(g)=a(\widetilde g)\}$. Define the constrained costs
\begin{equation}\label{eqsynchronized-distance}
\begin{split}
 W_{1,\mathrm{syn}}(\nu,\widetilde\nu)
 &=\inf_{\Gamma\in\operatorname{Cpl}_{\mathrm{syn}}(\nu,\widetilde\nu)}
                       \int d_\cG(g,\widetilde g)\dd\Gamma,\\
 \Wsyn(\boldsymbol\nu,\widetilde{\boldsymbol\nu})
 &=\int W_{1,\mathrm{syn}}(\nu_\omega,\widetilde\nu_\omega)\dd P(\omega).
\end{split}
\end{equation}
The second formula is used for kernels with identical translation
marginals almost everywhere. For the
state-dependent laws it is $\sum_s\pi_s W_{1,\mathrm{syn}}(\nu_s,
\widetilde\nu_s)$. On synchronized pairs, the cost is
$\min\{1,\|A(g)-A(\widetilde g)\|_\infty\}$.

\begin{lemma}[The constrained metric]\label[lemma]{lemsynchronized-metric}
On each class with fixed translation marginals, $\Wsyn$
defines a finite metric on kernels identified up to
$P$-almost-everywhere equality.
 One has
\begin{equation}\label{eqambient-vs-synchronized}
 \W\le\Wsyn\le1.
\end{equation}
There is in general no reverse comparison, even for laws
supported on a common compact set and sharing the same
nonatomic translation marginal.
\end{lemma}
\begin{proof}
Disintegrate both laws over their common translation marginal and
couple the conditional laws within each fiber. This gives a feasible
coupling with cost at most one. The space $\cG$ is Polish.
Couplings of fixed marginals are tight, since the mass outside a product
of compact sets is bounded by the sum of the two marginal tails.
The marginal conditions and the closed constraint $a(g)=a(\widetilde g)$
are preserved by weak limits. Prokhorov compactness
\cite[Theorem~23.2]{Kallenberg2021} and the bounded continuous cost
therefore give an optimal coupling, without compact
support. For convergent pairs of marginals, the same tightness argument
applied to nearly optimal couplings proves lower semicontinuity of the
value on the closed equal-marginal domain. Hence it is Borel, and the
integral in \eqref{eqsynchronized-distance} is measurable. The ambient infimum is taken over a larger set, proving the first
inequality and separation. Symmetry is immediate. For the triangle
inequality, glue two synchronized couplings over their common middle
marginal. All three translations then agree, and the metric triangle
inequality for $d_\cG$ bounds the $(1,3)$ cost by the sum of the other
two. Take infima and integrate.

To see the failure of a reverse comparison, take
$K=\TT\times\{A_0,A_1\}$ with distinct constant matrices, and let
$\nu=m\otimes\tfrac12(\delta_{A_0}+\delta_{A_1})$.
Divide the torus into $N$ equal intervals. Define $\nu_N$ to use
$A_0$ on the first half and $A_1$ on the second half of each interval,
with Haar translation marginal. By transporting within each interval
while keeping the matrix unchanged, $W_1(\nu_N,\nu)\le1/N$.
In a synchronized coupling, however, half the mass at each
translation must change its matrix. Therefore
\[
 W_{1,\mathrm{syn}}(\nu_N,\nu)
 =\tfrac12\min\{1,\|A_0-A_1\|\}>0.
\]
This separates the metrics, not Lyapunov regularity conclusions.
\end{proof}

\begin{remark}[Finite atomic marginals]\label[remark]{remmetric-comparison}
If each fixed $\mu_s$ has finite support and $h_*>0$ is the minimum
separation between distinct atoms (set $h_*=1$ when all supports are
singletons), then
\[
 \Wsyn\le(1+h_*^{-1})\W.
\]
Indeed, retain the translation-matched part of an ambient optimal
coupling. The residual marginals at each atom have the same mass and
can be recoupled within that fiber. The replaced mass is at most
$W_1/h_*$ and its new cost at most that mass. Integrate the resulting
statewise estimate. 
\end{remark}

For the exterior-square comparisons used below, set
\[
 V=\wedge^2\RR^r,\qquad Q_2=\TT^d\times\PP(V),
\]
and equip $V$ with its induced Euclidean norm.  For
$y=(\theta,\hat w)\in Q_2$, define
\begin{equation}\label{eqexterior-growth}
 H_n^{(2),\boldsymbol\nu}(\omega,y)
 =\EE_{\boldsymbol\nu,\omega}
   \log\frac{\|\wedge^2 A_n(\xi,\theta)w\|}{\|w\|}.
\end{equation}
The common matrix and inverse bound $B$ gives the corresponding
one-step bound $2M=2\log B$.

For $0<\beta\le1$ and $C_{\mathrm{ph}}<\infty$, define
\begin{equation}\label{eqholder-matrix-bounds}
 \mathcal K_{B,\beta,C_{\mathrm{ph}}}
 =\bigl\{g\in\cG_B:\
 \|A(g)(\theta)-A(g)(\theta')\|
 \le C_{\mathrm{ph}}
 d_{\TT^d}(\theta,\theta')^\beta
 \quad(\theta,\theta'\in\TT^d)\bigr\}.
\end{equation}
The matrix fields are uniformly bounded and equicontinuous, and
uniform limits preserve the matrix, inverse, and H\"older bounds.
Since the translation factor is compact, Arzel\`a--Ascoli yields a
compact generator class.
Set
\[
\mathcal E_{B,\beta,C_{\mathrm{ph}}}
=\mathcal E_B\cap
\Prob(\mathcal K_{B,\beta,C_{\mathrm{ph}}})^S.
\]

\begin{proposition}[Finite-time comparisons in two regimes]
\label[proposition]{propfinite-comparison}
For a pair of current-state kernels, set
\[
 E_n(\boldsymbol\nu,\widetilde{\boldsymbol\nu})
 =\max\left\{ |\mathcal L_n(\boldsymbol\nu)
                    -\mathcal L_n(\widetilde{\boldsymbol\nu})|,
 \frac1n\|H_n^{\boldsymbol\nu}-H_n^{\widetilde{\boldsymbol\nu}}\|_\infty
 \right\}.
\]
The following estimates hold for $D\le1$:
\begin{enumerate}
\item If both kernels belong to
$\mathcal C_{\boldsymbol\mu}(B)$, set
$D=\Wsyn$, $\kappa=1$, and $a=2\log B$.
\item If both kernels belong to
$\Prob(\mathcal K_{B,\beta,C_{\mathrm{ph}}})^S$, allowing the
translations to vary, set $D=\W$, $\kappa=\beta$, and
$a=2\log B+1$.
\end{enumerate}
In either case,
\begin{equation}\label{eqfinite-comparison}
 E_n(\boldsymbol\nu,\widetilde{\boldsymbol\nu})
 \le C e^{an}D^\kappa \quad(n\ge1).
\end{equation}
The constant depends only on $B,\pi_{\min}$ and, in (ii),
$\beta,C_{\mathrm{ph}}$.  In (i) one also has
\[
 \frac1n\|H_n^{\boldsymbol\nu}-H_n^{\widetilde{\boldsymbol\nu}}\|_{L^1(P;C(Q))}
 \le 2B^{2n}D,\qquad
 \frac1n\|H_n^{(2),\boldsymbol\nu}
        -H_n^{(2),\widetilde{\boldsymbol\nu}}\|_\infty
 \le\frac{4}{\pi_{\min}}B^{4n}D.
\]
\end{proposition}
\begin{proof}
For each state choose an $\eta$-optimal synchronized coupling in (i),
where $\eta>0$, and an optimal unconstrained coupling in (ii), where
we set $\eta=0$. Use their product along a fixed environment path.
The one-step expected cost is at most $\pi_{\min}^{-1}(D+\eta)$.
No selection over an infinite state set is involved. Write the evaluated
coupled matrices as $C_j=A(g_j)(\theta_j)$ and
$\widetilde C_j=A(\widetilde g_j)(\widetilde\theta_j)$, and let
$\delta_j=\|C_j-\widetilde C_j\|$.
For the $n$-step products $P_n,\widetilde P_n$, telescoping gives
\[
 P_n-\widetilde P_n
 =\sum_{j<n}C_{n-1}\cdots C_{j+1}
            (C_j-\widetilde C_j)\widetilde C_{j-1}\cdots\widetilde C_0,
 \qquad
 \|P_n-\widetilde P_n\|\le B^{n-1}\sum_{j<n}\delta_j.
\]
The inverse bounds imply
$\min\{\|P_n\|,\|\widetilde P_n\|,\|P_n v\|,
\|\widetilde P_n v\|\}\ge B^{-n}$ for every unit $v$.
Using $|\log u-\log v|\le|u-v|/\min(u,v)$, each logarithmic norm
difference is at most $B^{2n-1}\sum_{j<n}\delta_j$.

In (i), the coupled phase histories agree exactly. Since
$\|A-\widetilde A\|_\infty\le2B\min\{1,\|A-\widetilde A\|_\infty\}$,
$\EE\delta_j\le2B\pi_{\min}^{-1}(D+\eta)$.
Taking expectations, dividing by $n$, and sending $\eta\downarrow0$
gives $2\pi_{\min}^{-1}B^{2n}D$, uniformly in the initial state and
environment. The same comparison integrated over the Haar phase
bounds $|\mathcal L_n(\boldsymbol\nu)-\mathcal L_n(\widetilde{\boldsymbol\nu})|$.
For the integrated estimate, set $e_s=\int d_\cG\dd\Gamma_s$.
Stationarity gives
\[
 \int_\Omega\sum_{j<n}e_{\omega_j}\dd P
 =n\sum_s\pi_s e_s\le n(D+\eta).
\]
The same pathwise bound and $\eta\downarrow0$ give the sharper
$L^1(P;C(Q))$ estimate in the statement.

In (ii), put
$\chi_j=d_{\TT^d}(a(g_j),a(\widetilde g_j))$. Then
$\EE\chi_j\le\pi_{\min}^{-1}D$ and
$d_{\TT^d}(\theta_j,\widetilde\theta_j)\le\sum_{i<j}\chi_i$.
By \eqref{eqholder-matrix-bounds} and concavity of $t^\beta$,
\[
 \EE\delta_j
 \le2B\pi_{\min}^{-1}D
       +C_{\mathrm{ph}}(j\pi_{\min}^{-1}D)^\beta.
\]
The preceding product comparison yields
\[
 E_n\le B^{2n-1}
 \left(2B\pi_{\min}^{-1}D
        +C_{\mathrm{ph}}\pi_{\min}^{-\beta}n^\beta D^\beta\right).
\]
Now use $D\le D^\beta$ for $D\le1$ and $n^\beta\le e^n$.

For the additional exterior estimate in (i), use the isometric
antisymmetric-tensor realization of the induced Euclidean exterior norm.
The identity
\[
 C\otimes C-\widetilde C\otimes\widetilde C
 =(C-\widetilde C)\otimes C
       +\widetilde C\otimes(C-\widetilde C)
\]
gives
\[
 \|\wedge^2 C\|,\|(\wedge^2 C)^{-1}\|\le B^2,
 \qquad \|\wedge^2 C-\wedge^2\widetilde C\|
                       \le2B\|C-\widetilde C\|.
\]
Using the preceding telescoping estimate and logarithmic lower bound
with $B^2$ in place of $B$, we bound the exterior logarithmic difference by
$2B^{4n-1}\sum_{j<n}\delta_j$. Its expectation is at most
$4n\pi_{\min}^{-1}B^{4n}(D+\eta)$. Let $\eta\downarrow0$.
\end{proof}

\subsection{Probabilistic inputs and exterior-square deviations}
\label{subsecexterior}
For $\boldsymbol\nu_0\in\mathcal E_B$ and $\mu_s=a_*\nu_{0,s}$,
\cref{lemtranslation-ergodicity} and
\eqref{eqambient-vs-synchronized} give
\[
 \mathcal C_{\boldsymbol\mu}(B)\subset\mathcal E_B,
 \qquad \W\le\Wsyn.
\]
When the reference has SVU, $\Wsyn$-convergence implies
$\W$-convergence, so the qualitative continuity of $L_1$ at
$\boldsymbol\nu_0$ given by \cref{thmqualcontinuity} also holds in the
synchronized class. Under the same hypothesis, \cref{thmuniform-ldt}
provides an ambient relative neighborhood for each fixed tolerance.
Restricting its probability estimates to a sufficiently small synchronized
relative neighborhood preserves their prefactors, exponential rates,
and starting scales.

For the exterior estimates, take any reference
$\boldsymbol\nu_0\in\mathcal E_B$ and write
\[
 L_0=L_1(\boldsymbol\nu_0),\qquad
 \Lambda_0=L_1(\boldsymbol\nu_0)+L_2(\boldsymbol\nu_0),\qquad
 a_n^{(2)}(\omega)
 =\sup_{y\in Q_2}H_n^{(2),\boldsymbol\nu_0}(\omega,y).
\]
The identity $\|\wedge^2 A_n\|=s_1(A_n)s_2(A_n)$ identifies the top
exponent of the reference exterior cocycle with $\Lambda_0$.  Moreover,
$|H_n^{(2),\boldsymbol\nu}|\le2Mn$ and $|\Lambda_0|\le2M$.

\begin{lemma}[Upper conditional expected exterior growth]
\label[lemma]{lemexterior-expected}
For the reference kernel above,
\begin{equation}\label{eqexterior-expected-limit}
 \frac{a_n^{(2)}}n\longrightarrow\Lambda_0
 \quad P\text{-almost surely and in }L^1(P).
\end{equation}
\end{lemma}

\begin{proof}
For $r=2$, the exterior cocycle is the scalar identity and the assertion
is immediate.  For $r\ge3$, the map
$g\mapsto(a(g),\wedge^2A(g))$ is continuous and
$\det(\wedge^2A)=(\det A)^{r-1}=1$.  Its sampled mixed base is an
ergodic factor of the original one, and its top exponent is
$\Lambda_0$ by \eqref{eqexterior-exponent}.  Apply the maximal-envelope
part of \cref{thmexpected} to this cocycle.  
\end{proof}

\begin{proposition}[Exterior upper deviations uniform over nearby kernels]
\label[proposition]{propexterior-ldt}
Let $\boldsymbol\nu_0$ be as above and consider either of the following
perturbation settings:
\begin{enumerate}
\item the synchronized class $\mathcal C_{\boldsymbol\mu}(B)$, where
      $\mu_s=a_*\nu_{0,s}$ for every $s$, with metric $\Wsyn$;
\item the bounded ambient class $\mathcal E_B$, where
      $\boldsymbol\nu_0\in\mathcal E_B$, with metric $\W$.
\end{enumerate}
For every $\tau>0$ there exist a relative neighborhood
$\mathcal V_\tau$ of $\boldsymbol\nu_0$ in the corresponding class,
constants $C,c>0$, and $n_0\ge1$ such that
\begin{equation}\label{eqexterior-upper-ldt}
\sup_{\boldsymbol\nu\in\mathcal V_\tau}
\int_\Omega\sup_{\theta\in\TT^d}
\PP_{\boldsymbol\nu,\omega}
\left\{\frac1n\log\|\wedge^2 A_n(\xi,\theta)\|
                         >\Lambda_0+\tau\right\}\,dP(\omega)
\le Ce^{-cn}
\qquad(n\ge n_0).
\end{equation}
\end{proposition}

\begin{proof}
For $r=2$, the deviation events are empty.  Suppose $r\ge3$, put
$\varepsilon=\tau/2$, and define
\[
 u_\ell^{\boldsymbol\nu}(\omega)
 =\left(\frac1\ell\sup_{y\in Q_2}
 H_\ell^{(2),\boldsymbol\nu}(\omega,y)-\Lambda_0\right)_+.
\]
By \cref{lemexterior-expected}, choose $\ell$ such that
$\int u_\ell^{\boldsymbol\nu_0}\,dP\le\varepsilon/16$.
Exterior fixed-time continuity follows from
\cref{lemfinite-continuity}, with restriction to synchronized
neighborhoods as above. Since the
positive-part map is $1$-Lipschitz, in either setting we may shrink a
relative neighborhood so that
\[
 \|u_\ell^{\boldsymbol\nu}
       -u_\ell^{\boldsymbol\nu_0}\|_\infty
 \le\frac1\ell
 \|H_\ell^{(2),\boldsymbol\nu}
       -H_\ell^{(2),\boldsymbol\nu_0}\|_\infty
 \le\frac{\varepsilon}{16}.
\]
Hence $\int u_\ell^{\boldsymbol\nu}\,dP\le\varepsilon/8$ throughout
the neighborhood.  The center remains the fixed number $\Lambda_0$.

Each $u_\ell^{\boldsymbol\nu}$ is a length-$\ell$ word observable with
values in $[0,4M]$.  For
\[
 E_n^{\boldsymbol\nu}
 =\left\{\omega:
 \frac1n\sum_{j=0}^{n-1}
 u_\ell^{\boldsymbol\nu}(\Theta^j\omega)>\frac\varepsilon4\right\},
\]
\cref{lemword-concentration}, with range bound $4M$ and deviation
$\varepsilon/8$, gives uniformly in the chosen neighborhood
\[
 P(E_n^{\boldsymbol\nu})
 \le\exp\left(-\frac{\varepsilon^2n}
                    {32768M^2T_\ell^2}\right)
 \quad\text{for }n\ge
 \left\lceil\frac{128MT_\ell}{\varepsilon}\right\rceil.
\]
Apply \eqref{eqblock-one-sided} with one-step bound $2M$, center
$\Lambda_0$, and $t=\varepsilon/2$.  For
$n\ge\max\{\ell,\lceil16M\ell/\varepsilon\rceil\}$ the boundary term
is at most $\varepsilon/4$, so on $(E_n^{\boldsymbol\nu})^c$
\[
\sup_{y=(\theta,\hat w)\in Q_2}
\PP_{\boldsymbol\nu,\omega}
\left\{\frac1n\log
 \frac{\|\wedge^2 A_n(\xi,\theta)w\|}{\|w\|}
 >\Lambda_0+\varepsilon\right\}
\le\ell\exp\left(-\frac{\varepsilon^2n}{128M^2\ell}\right).
\]
Bound this probability by $1$ on $E_n^{\boldsymbol\nu}$ and integrate.

Finally, let $w_1,\ldots,w_N$ be an orthonormal basis of $V$, where
$N=\dim V$.  Since $\|T\|\le\sqrt N\max_i\|Tw_i\|$, once
$\log N/(2n)\le\tau/2$ the operator event in
\eqref{eqexterior-upper-ldt} is contained, for each fixed phase, in the
union of the $N$ directional events at threshold
$\Lambda_0+\varepsilon$.  Taking probabilities, the supremum over the
phase, and then the environmental integral proves the proposition after
adjusting $C,c,n_0$.
\end{proof}

\subsection{Uniform gaps and corrected scales}
\label{subseccommon-gap}
\label{subseccorrected-scales}

\begin{proposition}[Common finite-scale gap estimates]
\label[proposition]{propcommon-gap}
Retain $L_0$ and $\Lambda_0$ from the reference kernel above, and work
in either perturbation setting of \cref{propexterior-ldt}.  Suppose
that the reference kernel $\boldsymbol\nu_0$ has \SVU\ and
\begin{equation}\label{eqreference-gap}
 g_0=L_1(\boldsymbol\nu_0)-L_2(\boldsymbol\nu_0)>0.
\end{equation}
For every fixed $0<\tau<g_0/100$, there exist a relative neighborhood
$\mathcal U$ of $\boldsymbol\nu_0$ in the corresponding parameter
class, constants $C_0,c_0>0$, and $n_0\ge1$ such that, for every
$\boldsymbol\nu\in\mathcal U$ and $n\ge n_0$,
\begin{align}
\rho_{\boldsymbol\nu}
\left\{\left|\frac1n\log\|A_n\|-L_0\right|>\tau\right\}
&\le C_0e^{-c_0n},\label{eqcommon-top-ldt}\\
\rho_{\boldsymbol\nu}
\left\{\frac1n\log\|\wedge^2 A_n\|>\Lambda_0+\tau\right\}
&\le C_0e^{-c_0n}.\label{eqcommon-wedge-ldt}
\end{align}
Moreover, every kernel in this neighborhood satisfies
\begin{equation}\label{eqnearby-gap}
 L_1(\boldsymbol\nu)-L_2(\boldsymbol\nu)
 \ge g_0-3\tau>0.
\end{equation}
\end{proposition}

\begin{proof}
Apply the operator-norm conclusion of \cref{thmuniform-ldt} at
$\tau/2$, restricting to synchronized neighborhoods as in
Section~\ref{subsecexterior} when appropriate. By
\cref{thmqualcontinuity}, shrink the neighborhood so that
\[
 |L_1(\boldsymbol\nu)-L_0|<\frac\tau2.
\]
Then
\[
\left\{\left|\frac1n\log\|A_n\|-L_0\right|>\tau\right\}
\subset
\left\{\left|\frac1n\log\|A_n\|-L_1(\boldsymbol\nu)\right|
 >\frac\tau2\right\},
\]
and integration over the Haar phase gives \eqref{eqcommon-top-ldt}.
Intersect this neighborhood with the corresponding
$\mathcal V_\tau$ from \cref{propexterior-ldt} to obtain
\eqref{eqcommon-wedge-ldt}; take the larger prefactor and starting scale
and the smaller decay rate.

Set $W_n=n^{-1}\log\|\wedge^2A_n\|$ and
$D_n=\{W_n>\Lambda_0+\tau\}$.  Since
$|W_n|,|\Lambda_0|\le2M$,
\[
 W_n\le\Lambda_0+\tau+(4M+\tau)\mathbf1_{D_n},
 \qquad
 \int W_n\,d\rho_{\boldsymbol\nu}
 \le\Lambda_0+\tau+(4M+\tau)C_0e^{-c_0n}.
\]
Letting $n\to\infty$ and using \eqref{eqexterior-exponent} gives
\[
 L_1(\boldsymbol\nu)+L_2(\boldsymbol\nu)
 \le\Lambda_0+\tau.
\]
The chosen neighborhood also gives
$L_1(\boldsymbol\nu)\ge L_0-\tau$, and hence
\[
\begin{aligned}
 L_1(\boldsymbol\nu)-L_2(\boldsymbol\nu)
 &=2L_1(\boldsymbol\nu)
   -[L_1(\boldsymbol\nu)+L_2(\boldsymbol\nu)]\\
 &\ge2(L_0-\tau)-(\Lambda_0+\tau)
   =g_0-3\tau>0.
\end{aligned}
\]
\end{proof}

The deterministic input is the Avalanche Principle
\cite[Proposition~4.2 in the chapter preprint]{DKbook}.

There are constants $c_{\mathrm{AP}},C_{\mathrm{AP}}>0$ such that,
for $p\ge3$, $0<\epsilon<1$, and
$0<\varkappa\le c_{\mathrm{AP}}\epsilon^2$, matrices
$B_0,\ldots,B_{p-1}\in\mathrm{GL}_r(\RR)$ satisfying
\begin{equation}\label{eqap-hypotheses}
\begin{aligned}
\frac{s_1(B_j)}{s_2(B_j)}&>\varkappa^{-1}
&& (0\le j<p),\\
\frac{\|B_{j+1}B_j\|}{\|B_{j+1}\|\,\|B_j\|}&>\epsilon
&& (0\le j<p-1)
\end{aligned}
\end{equation}
obey
\begin{equation}\label{eqap-statement}
\begin{split}
\left|\log\|B_{p-1}\cdots B_0\|
 +\sum_{j=1}^{p-2}\log\|B_j\|
 -\sum_{j=0}^{p-2}\log\|B_{j+1}B_j\|\right|
\le C_{\mathrm{AP}}p\frac{\varkappa}{\epsilon^2}.
\end{split}
\end{equation}
Here $s_1\ge s_2$ are the two largest singular values. Define
\begin{equation}\label{eqcorrected-scale}
\mathcal F_n(\boldsymbol\nu)
=2\mathcal L_{2n}(\boldsymbol\nu)-\mathcal L_n(\boldsymbol\nu).
\end{equation}
The exponential accuracy below concerns these corrected averages,
not the ordinary averages $\mathcal L_n$.

\begin{proposition}[Exponentially accurate corrected scales]
\label[proposition]{propcorrected-scales}
Under the hypotheses of \cref{propcommon-gap}, fix $\tau$ and the
neighborhood $\mathcal U$ furnished there. After increasing $n_0$,
there are constants $C,b>0$ such that
\begin{equation}\label{eqcorrected-approximation}
\sup_{\boldsymbol\nu\in\mathcal U}
|L_1(\boldsymbol\nu)-\mathcal F_n(\boldsymbol\nu)|
\le Ce^{-bn}\qquad(n\ge n_0).
\end{equation}
\end{proposition}

\begin{proof}
Fix $\boldsymbol\nu\in\mathcal U$ and abbreviate
$\rho=\rho_{\boldsymbol\nu}$,
$\mathcal L_k=\mathcal L_k(\boldsymbol\nu)$, and
$\mathcal F_k=\mathcal F_k(\boldsymbol\nu)$.
\emph{Good blocks and the AP hypotheses.}
For $n\ge n_0$ and $p\ge3$, put
$B_j(z)=A_n(F^{jn}z)$ for $0\le j<p$. Then
$B_{p-1}\cdots B_0=A_{pn}(z)$ and
$B_{j+1}B_j=A_{2n}(F^{jn}z)$.
Let $G_{n,p}$ be the event on which
\[
\begin{aligned}
\left|\frac1n\log\|B_j\|-L_0\right|&\le\tau,
&\frac1n\log\|\wedge^2 B_j\|&\le\Lambda_0+\tau
&& (0\le j<p),\\
\left|\frac1{2n}\log\|B_{j+1}B_j\|-L_0\right|&\le\tau
&&&& (0\le j<p-1).
\end{aligned}
\]
Invariance and \cref{eqcommon-top-ldt,eqcommon-wedge-ldt} give
\begin{equation}\label{eqap-bad-event}
\begin{aligned}
\rho(G_{n,p}^c)
&\le2pC_0e^{-c_0n}+(p-1)C_0e^{-2c_0n}\\
&\le3C_0p e^{-c_0n}.
\end{aligned}
\end{equation}
This uses neither block independence nor ergodicity of $F^n$.
On $G_{n,p}$,
\begin{equation}\label{eqap-gap-calculation}
\begin{aligned}
\log\frac{s_1(B_j)}{s_2(B_j)}
&=2\log\|B_j\|-\log\|\wedge^2 B_j\|\\
&\ge2n(L_0-\tau)-n(\Lambda_0+\tau)
 =n(g_0-3\tau).
\end{aligned}
\end{equation}
The single- and double-block estimates also give
\begin{equation}\label{eqap-noncancellation}
\log\frac{\|B_{j+1}B_j\|}{\|B_{j+1}\|\,\|B_j\|}
\ge2n(L_0-\tau)-2n(L_0+\tau)=-4\tau n.
\end{equation}
Take
\begin{equation}\label{eqap-parameters}
\kappa_n=e^{-(g_0-4\tau)n},
\qquad
\epsilon_n=e^{-5\tau n},
\qquad
a_0=g_0-14\tau>0.
\end{equation}
The gap ratio exceeds $\kappa_n^{-1}$, the adjacent-block ratio exceeds
$\epsilon_n$, and
$\kappa_n/\epsilon_n^2=e^{-a_0n}$.
Increase $n_0$ so that this last quantity is at most $c_{\mathrm{AP}}$.
Thus \cref{eqap-hypotheses} holds with
$\varkappa=\kappa_n$ and $\epsilon=\epsilon_n$.

\emph{The integrated estimate.}
Write the signed AP expression as
\[
\mathcal R_{n,p}
=\log\|A_{pn}\|
 +\sum_{j=1}^{p-2}\log\|B_j\|
 -\sum_{j=0}^{p-2}\log\|B_{j+1}B_j\|.
\]
On $G_{n,p}$, the AP gives
$|\mathcal R_{n,p}|\le C_{\mathrm{AP}}p e^{-a_0n}$.
Everywhere, the common norm bounds give
\[
|\mathcal R_{n,p}|
\le Mpn+(p-2)Mn+2(p-1)Mn\le4Mpn.
\]
Splitting the integral over the good and bad events yields
\[
\frac1{pn}\left|\int\mathcal R_{n,p}\,d\rho\right|
\le\frac{C_{\mathrm{AP}}}{n}e^{-a_0n}+4M\rho(G_{n,p}^c)
\le C\bigl(e^{-a_0n}+p e^{-c_0n}\bigr).
\]
By invariance, the single-, double-, and $p$-block terms integrate to
$n\mathcal L_n$, $2n\mathcal L_{2n}$, and $pn\mathcal L_{pn}$,
respectively. Hence
\begin{equation}\label{eqap-integrated}
\left|\mathcal L_{pn}
 +\left(1-\frac2p\right)\mathcal L_n
 -2\left(1-\frac1p\right)\mathcal L_{2n}\right|
\le C\bigl(e^{-a_0n}+p e^{-c_0n}\bigr).
\end{equation}
The expression inside the absolute value is
$\mathcal L_{pn}-\mathcal F_n
 +(2/p)(\mathcal L_{2n}-\mathcal L_n)$.
Since $|\mathcal L_k|\le M$, the correction is at most $4M/p$, so
\begin{equation}\label{eqap-finite-relation}
|\mathcal L_{pn}-\mathcal F_n|
\le C\bigl(p^{-1}+e^{-a_0n}+p e^{-c_0n}\bigr),
\end{equation}
with constants independent of the kernel, $p$, and $n$.

\emph{Scale iteration.}
Set $b=\tfrac14\min\{a_0,c_0\}>0$ and increase $n_0$ so that
$p(n)=\lfloor e^{bn}\rfloor\ge3$. Then
\[
p(n)^{-1}\le2e^{-bn},\qquad
e^{-a_0n}\le e^{-4bn},\qquad
p(n)e^{-c_0n}\le e^{-3bn}.
\]
Applying \cref{eqap-finite-relation} with $p(n)$ and $2p(n)$, and
writing $N=np(n)$, gives
\[
|\mathcal L_N-\mathcal F_n|
 +|\mathcal L_{2N}-\mathcal F_n|\le Ce^{-bn},
\]
and therefore
\[
|\mathcal F_N-\mathcal F_n|
\le2|\mathcal L_{2N}-\mathcal F_n|
    +|\mathcal L_N-\mathcal F_n|
\le Ce^{-bn}.
\]
Starting at an arbitrary $n\ge n_0$, define
\[
n_0'=n,\qquad n_{j+1}'=n_j'\lfloor e^{bn_j'}\rfloor.
\]
Since $n_j'\ge3^j n$, iteration gives
\begin{equation}\label{eqcorrected-telescoping}
|\mathcal F_{n_k'}-\mathcal F_n|
\le C\sum_{j=0}^{k-1}e^{-bn_j'}
\le C'e^{-bn}.
\end{equation}
Indeed, the sum is at most
$e^{-bn}/(1-e^{-2bn_0})$, using $3^j\ge1+2j$.
For each fixed kernel, \cref{eqlyapunov} gives
$\mathcal F_{n_k'}(\boldsymbol\nu)\to L_1(\boldsymbol\nu)$.
Since the telescoping constants are kernel-independent, passing to the
limit proves \cref{eqcorrected-approximation} uniformly without assuming
prior uniform convergence.
\end{proof}

\subsection{H\"older continuity}
\label{subsecregularity}
The local two-point H\"older bound for $L_1$ follows by balancing the
corrected-scale approximation error and the finite-time kernel
comparison error in \eqref{eqholder-balance}.

\begin{theorem}[Local two-point H\"older continuity]\label[theorem]{thmholder}
Fix the finite-state irreducible aperiodic Markov environment. Let
$\boldsymbol\nu_0\in\mathcal E_B$ have \SVU\ and
$g_0=L_1(\boldsymbol\nu_0)-L_2(\boldsymbol\nu_0)>0$.
The following conclusions hold.
\begin{enumerate}
\item Let $\mu_s=a_*\nu_{0,s}$. On
$\mathcal C_{\boldsymbol\mu}(B)$ there are $\delta,C,\gamma>0$ such
that
\begin{equation}\label{eqholder-fixed}
 |L_1(\boldsymbol\nu)-L_1(\widetilde{\boldsymbol\nu})|
 \le C\Wsyn(\boldsymbol\nu,\widetilde{\boldsymbol\nu})^\gamma
\end{equation}
whenever both kernels have constrained distance less than $\delta$
from $\boldsymbol\nu_0$.
\item Fix $0<\beta\le1$ and $C_{\mathrm{ph}}<\infty$, and assume
$\boldsymbol\nu_0\in
\mathcal E_{B,\beta,C_{\mathrm{ph}}}$. There are
$\delta,C,\gamma>0$ such that
\begin{equation}\label{eqholder-varying}
 |L_1(\boldsymbol\nu)-L_1(\widetilde{\boldsymbol\nu})|
 \le C\W(\boldsymbol\nu,\widetilde{\boldsymbol\nu})^\gamma
\end{equation}
for kernels in $\mathcal E_{B,\beta,C_{\mathrm{ph}}}$ at ambient
distance less than $\delta$ from $\boldsymbol\nu_0$.
\end{enumerate}
The constants depend on the reference kernel, its gap, the fixed chain,
the dimensions, and $B$; in (ii) they also depend on
$\beta,C_{\mathrm{ph}}$. 
\end{theorem}
\begin{proof}
Choose $0<\tau<g_0/100$. In each regime,
\cref{propcommon-gap,propcorrected-scales} give a relative neighborhood
of the reference on which
\[
 |L_1(\boldsymbol\nu)-\mathcal F_n(\boldsymbol\nu)|
 \le C_0e^{-bn}\quad(n\ge n_0).
\]
In (i), choose a relative $\Wsyn$ ball in this neighborhood. In (ii),
intersect the ambient neighborhood with
$\mathcal E_{B,\beta,C_{\mathrm{ph}}}$ and choose a relative $\W$ ball. Take both radii at most $1/2$ and write
$(D,\kappa)=(\Wsyn,1)$ or $(\W,\beta)$, respectively. In either case $D\le1$ for pairs
in the chosen ball.

For the corresponding $a>0$, \cref{propfinite-comparison} gives
\[
 |\mathcal F_n(\boldsymbol\nu)
             -\mathcal F_n(\widetilde{\boldsymbol\nu})|
 \le 3C_1e^{2an}D^\kappa.
\]
Putting $a_F=2a$, the three-term comparison with the two limiting
exponents gives
\begin{equation}\label{eqholder-balance}
 |L_1(\boldsymbol\nu)-L_1(\widetilde{\boldsymbol\nu})|
 \le2C_0e^{-bn}+3C_1e^{a_Fn}D^\kappa\quad(n\ge n_0).
\end{equation}
For $0<D\le1$ take
\[
 n=\left\lfloor\frac{\kappa}{a_F+b}\log\frac1D\right\rfloor,
 \qquad \gamma=\frac{\kappa b}{a_F+b}>0.
\]
If $n\ge n_0$, then $e^{-bn}\le e^b D^\gamma$ and
$e^{a_F n}D^\kappa\le D^\gamma$. If $n<n_0$, then
$D>\exp(-(a_F+b)n_0/\kappa)$ and the bound
$|L_1(\boldsymbol\nu)-L_1(\widetilde{\boldsymbol\nu})|\le2M$
is absorbed by increasing $C$. The case $D=0$ is immediate.
All constants are independent of the two kernels in the ball.
\end{proof}

\begin{remark}[Scope]\label[remark]{remholder-dependence}
Because the reference block scale is chosen
through qualitative expected-growth convergence, the H\"older exponent
is not asserted to be optimal or effectively computable from SVU alone.
For fixed finite atomic translation marginals, \cref{remmetric-comparison}
also converts (i) to ambient Wasserstein H\"older continuity, after reducing
the neighborhood and changing the constant. 
\end{remark}

\Needspace{7\baselineskip}
\begin{example}[Fixed irrational phase with non-H\"older matrix fields]
\label{ex:noncoverage}
Take $d=1$, $r=2$, a fixed positive two-state transition matrix
$\Pi$, and an irrational $\alpha\in\TT$. Let $f\ge0$ be a
continuous circle function, supported near $0$, with $f(0)=0$ and
$f(t)=1/\log(e/|t|)$ for sufficiently small $t\ne0$.
A smooth cutoff of this local formula gives such a function. Fix distinct
$c_1,c_2>0$, put $a_s(\theta)=c_s+f(\theta)$, and let $\nu_s$
be the image of Haar measure on $\mathrm{SO}(2)$ under
\[
 R\longmapsto g_{s,R}=(\alpha,A_{s,R}),\qquad
 A_{s,R}(\theta)=
 \operatorname{diag}(e^{a_s(\theta)},e^{-a_s(\theta)})R.
\]
These continuous images of $\mathrm{SO}(2)$ are compact in $\cG$,
with matrix and inverse norms bounded by
$\exp(\max_s c_s+\|f\|_\infty)$. The environmental shift is mixing,
so its product with the irrational rotation is ergodic: nonzero
phase Fourier coefficients of an invariant function would be
nontrivial eigenfunctions of the mixing shift. By the harmonic-function
criterion in \cref{lemtranslation-ergodicity}, the mixed base is
therefore ergodic.

Haar rotation invariance and \eqref{eqrenewal-one-step} give
$\int\Psi(g,(\theta,\hat p))\dd\nu_s(g)=\log\cosh a_s(\theta)$,
independently of $\hat p$. Every stationary projective family has
Haar phase marginal by \cref{thmstationary}, so all stationary values equal
\[
 L_1=\sum_s\pi_s\int_\TT\log\cosh a_s(\theta)\dd m(\theta)>0.
\]
Thus SVU holds and $L_2=-L_1<L_1$. On the other hand, right
multiplication by $R$ preserves the operator norm, so
\[
 \|A_{s,R}(t)-A_{s,R}(0)\|
 \asymp f(t)=\frac1{\log(e/|t|)}\qquad(t\to0).
\]
Each matrix field is continuous but not H\"older of any positive
exponent in the usual phase metric. Thus \cref{thmholder}(i) applies
in the bounded synchronized class with fixed marginals
$(\delta_\alpha,\delta_\alpha)$. That branch does not require these supports to be compact.

In the natural enlarged Markov encoding the phase law is
$\delta_{\theta_0+n\alpha}$, whose total-variation distance from
Haar measure is $1$. For a nonzero Fourier mode $e_k$, the phase
operator has $\|Q^n e_k\|_\infty=1$. Thus the uniform-ergodicity condition discussed in \cite[Remark~2.1]{CDKM} and the strong-mixing estimate (1.4) required in Definition~1.1 of \cite{CDKstatistics} do not follow from these hypotheses: for every nonzero Fourier mode, $\|Q^n e_k\|_\infty=1$.
This tests direct applicability to that encoding, not every possible
reduction or observable-specific use of those results.
\end{example}

\section{Moving-law Schr\"odinger cocycles}\label{secschrodinger}
Moving-law Schr\"odinger cocycles with continuous potentials have
H\"older energy dependence under the same-frequency hypotheses below. Energy changes preserve sampled
translations, so the synchronized theorem applies without phase
H\"older regularity. Positivity and strong irreducibility use the
external criteria identified in the proof, whereas mixed-base
ergodicity is verified directly.

Let $v\in C^0(\TT^d,\RR)$, $R<\infty$, and for each state $s$ fix
$\kappa_s\in\Prob(\TT^d\times[-R,R])$. Frequency and noise may be
correlated. Put
\begin{equation}\label{eqschrodinger-matrix}
 A_{E,\zeta}(\theta)=
 \begin{pmatrix}v(\theta)+\zeta-E&-1\\1&0\end{pmatrix},
 \qquad
 \nu_{s,E}=\bigl((\alpha,\zeta)\mapsto(\alpha,A_{E,\zeta})\bigr)_*
                                                        \kappa_s,
\end{equation}
and write $\boldsymbol\nu_E=(\nu_{s,E})_{s\in S}$ and
$\mu_s=(\operatorname{pr}_{\TT^d})_*\kappa_s$.

\begin{theorem}[H\"older energy dependence]\label[theorem]{thmschrodinger}
Let the environment be the stationary two-sided shift of a fixed
finite irreducible aperiodic transition matrix $\Pi$, with stationary
law $\pi$. Suppose:
\begin{enumerate}
\item some $\alpha_*\in\bigcap_s\supp\mu_s$ satisfies
$e^{2\pi i k\cdot\alpha_*}\ne1$ for every
$k\in\ZZ^d\setminus\{0\}$;
\item there is $S_{\mathrm{sf}}\subset S$ such that each
$s\in S_{\mathrm{sf}}$ has two support points
$(\alpha_s,\zeta_{1,s}), (\alpha_s,\zeta_{2,s})\in\supp\kappa_s$
with $\zeta_{1,s}\ne\zeta_{2,s}$, and
\begin{equation}\label{eqadjacent-noisy-states}
 \sum_{s,t\in S_{\mathrm{sf}}}\pi_s\Pi_{st}>0.
\end{equation}
\end{enumerate}
Then the mixed base is ergodic and $L_1(\boldsymbol\nu_E)>0$ for every
real $E$. For each compact energy interval $I$ there are
$C_I<\infty$ and $\gamma_I>0$ such that
\begin{equation}\label{eqenergy-holder}
 |L_1(\boldsymbol\nu_E)-L_1(\boldsymbol\nu_{E'})|
 \le C_I|E-E'|^{\gamma_I}\quad(E,E'\in I).
\end{equation}
In particular $\inf_{E\in I}L_1(\boldsymbol\nu_E)>0$.
\end{theorem}
\begin{proof}
The finite-state Markov shift is mixing, hence weakly mixing.
The direct verification of mixed-base ergodicity below uses
\cref{lemtranslation-ergodicity} and specializes
\cite[Proposition~6.2]{LCZpositivity} to the finite-state setting. If a bounded real function $h$ is
harmonic for the environment--phase operator $K$, Jensen's equality
gives transport invariance along its transitions. Its Fourier
coefficients $c_k(\omega)$ consequently satisfy
\[
 c_k(\omega)=e^{2\pi i k\cdot\alpha}
                         c_k(\Theta\omega)
 \quad\text{for }\mu_{\omega_0}\text{-almost every }\alpha.
\]
Continuity in $\alpha$ extends this identity to the support,
including $\alpha_*$. Thus $c_k\circ\Theta
=e^{-2\pi i k\cdot\alpha_*}c_k$. Weak mixing and (i) force
$c_k=0$ for $k\ne0$, while $c_0$ is constant by ergodicity.
Fourier uniqueness makes $h$ constant, proving the claim.
This argument is independent of the energy.

Let $\Omega_{\mathrm{sf}}=\{\omega:\omega_0\in S_{\mathrm{sf}}\}$.
The adjacent noisy-environment overlap is
\[
 P\{\omega_{-1},\omega_0\in S_{\mathrm{sf}}\}
 =\sum_{s,t\in S_{\mathrm{sf}}}\pi_s\Pi_{st}>0.
\]
Together with (i), bounded noise, and continuity of $v$, this overlap
allows the same-frequency criterion
\cite[Proposition~6.4]{LCZpositivity} to give
$L_1(\boldsymbol\nu_E)>0$ at every energy. Its monodromy verification
and \cite[Theorem~5.11]{LCZpositivity} give environment-wise strong
irreducibility. More explicitly, distinct noises at a current
state and at its preceding state produce respectively a nontrivial
lower and upper unipotent return at the same environment--phase
fiber. The cited exclusion theorem converts these matched returns
into strong irreducibility. By \cref{corstrong-SVU}, every energy kernel has \SVU. Since $\det A_{E,\zeta}=1$,
$L_2(\boldsymbol\nu_E)=-L_1(\boldsymbol\nu_E)<L_1(\boldsymbol\nu_E)$.

Fix a compact $J$ with $I\subset\operatorname{int}J$ and set
\[
 K_J=\{(\alpha,A_{E,\zeta}):E\in J,
                 (\alpha,\zeta)\in\bigcup_s\supp\kappa_s\}.
\]
This compact image has a common matrix and inverse bound $B_J\ge e$.
Thus $\boldsymbol\nu_E\in\mathcal C_{\boldsymbol\mu}(B_J)$ and
$|L_1(\boldsymbol\nu_E)|\le\log B_J$ for $E\in J$.
Couple $\nu_{s,E}$ and $\nu_{s,E'}$ by using the same
$(\alpha,\zeta)$. Since
\[
 A_{E,\zeta}-A_{E',\zeta}
 =\begin{pmatrix}E'-E&0\\0&0\end{pmatrix},
\]
this is a synchronized coupling with cost at most $|E-E'|$.
Therefore
\begin{equation}\label{eqenergy-wasserstein}
 \Wsyn(\boldsymbol\nu_E,\boldsymbol\nu_{E'})\le|E-E'|.
\end{equation}
The first regime of \cref{thmholder} now gives a two-point H\"older
estimate on an open interval about every $E_0\in I$.

Choose a finite subcover $U_1,\ldots,U_N$ with local constants
$C_j,\gamma_j$, and a Lebesgue number $0<\delta_0\le1$.
Set $\gamma_I=\min_j\gamma_j>0$. Pairs with $|E-E'|<\delta_0$
lie in one $U_j$; for the other pairs use the bound $2\log B_J$.
Then \eqref{eqenergy-holder} holds with
\[
 C_I=\max\{\max_j C_j,\,2\log B_J\,\delta_0^{-\gamma_I}\}.
\]
Strict positivity and continuity give a positive minimum on $I$.
\end{proof}

\begin{remark}[Weaker environments]
For a general ergodic environment and bounded noise, the same
energy coupling and common matrix bound give continuity at every
\SVU\ energy by \cref{thmqualcontinuity}, without a
finite-state hypothesis. However, this qualitative argument does not provide a H\"older modulus.
\end{remark}

\begin{example}[A concrete example]
Take
$\Pi=\left(\begin{smallmatrix}3/4&1/4\\1/3&2/3\end{smallmatrix}\right)$,
$\pi=(4/7,3/7)$, a rationally independent $\alpha_*$, and a distinct
$\alpha^\dagger$. The laws
\[
 \kappa_1=\tfrac13\delta_{(\alpha_*,-1)}
          +\tfrac13\delta_{(\alpha_*,1)}
          +\tfrac13\delta_{(\alpha^\dagger,0)},\qquad
 \kappa_2=\tfrac12\delta_{(\alpha_*,2)}
          +\tfrac12\delta_{(\alpha^\dagger,-2)}
\]
satisfy the theorem with $S_{\mathrm{sf}}=\{1\}$, since
$\pi_1\Pi_{11}=3/7>0$. 
\end{example}

\appendix
\section{Measurable subbundles and support exactification}
Measurable subspace fields admit Grassmannian representations and
measurable frames. Their almost-sure invariance extends to every
generator in the conditional support, with the precise exceptional-set
quantifiers needed for \cref{corstrong-SVU}.

\subsection{Subbundle coding and measurable frames}
\label{app:subbundle-coding}

Grassmannian maps and measurable total spaces give equivalent
representations of the subbundles in
\cref{rem:equivalent-subbundle-descriptions}. For a measurable
$L:\Omega\times\TT^d\to\Gr_k(\RR^r)$, the total space is measurable
because the tautological set
\[
 \mathcal T_k=\{(W,v)\in\Gr_k(\RR^r)\times\RR^r:v\in W\}
\]
is closed.

Conversely, suppose a subspace field $L$ has a measurable graph with
closed linear sections, as does $L_c$ in \cref{lem:slow-subspaces}.
For a countable dense set $(q_j)$ in $\RR^r$ and $t>0$,
\[
\begin{aligned}
 &\{(\omega,\theta):\dist(q_j,L(\omega,\theta))<t\}\\
 &\quad=\bigcup_{\substack{s\in\mathbb Q\\0<s<t}}
 \pi_{\Omega\times\TT^d}
 \left(\operatorname{graph}(L)\cap
 ((\Omega\times\TT^d)\times\overline B(q_j,s))\right).
\end{aligned}
\]
The sets under projection are Borel with compact vertical sections, so
Arsenin--Kunugui \cite[Theorem~5.12.1]{Srivastava1998} makes each
projection Borel. These measurable distance functions determine the
closed subspace and hence its point in
$\bigsqcup_{\ell=0}^r\Gr_\ell(\RR^r)$. Thus the field is measurable,
including its rank. If that rank equals $k$ almost everywhere,
assigning a fixed $W_*\in\Gr_k(\RR^r)$ on the measurable exceptional
set gives a pointwise fixed-rank version.

On each positive-rank stratum, the continuous projection from the
compact Stiefel space onto the Grassmannian admits a measurable section
by \cite[Corollary~5.2.6]{Srivastava1998}. This supplies the measurable
orthonormal frames used in \cref{lem:slow-subspaces} and the
incidence-space construction of \cref{lem:bundle-stationary-family}.

\subsection{Support exactification}\label{appsupport}
Almost-sure subbundle invariance extends to every generator in the
conditional support, while the exceptional phase set may depend on
$(\omega,g)$. Continuity of the action on Dirac Young measures gives
this extension without assuming continuity of the subbundle itself.

For $1\leq k<r$, set $Q_k:=\TT^d\times\Gr_k(\RR^r)$ and define
\begin{equation}\label{eq:grassmann-action}
  g\cdot(\theta,W)
  :=(\theta+a(g),A(g)(\theta)W).
\end{equation}
This is a continuous action of $\cG$ on the compact metric space $Q_k$.

For a measurable $k$-dimensional subbundle $L$, define
\begin{equation}\label{eq:dirac-young-bundle}
 \Lambda_\omega^L(\dd\theta,\dd W)
 :=m(\dd\theta)\delta_{L(\omega,\theta)}(\dd W).
\end{equation}
The map $\omega\mapsto\Lambda_\omega^L$ is a measurable probability kernel on $Q_k$.

\begin{lemma}[Grassmannian support exactification]\label[lemma]{lem:grassmannian-support-upgrade}
Let $L$ be a $k$-dimensional relative invariant subbundle.  Then there exists a Borel set $\Omega_L\subset\Omega$ with $P(\Omega_L)=1$ such that, for every $\omega\in\Omega_L$ and every $g\in\supp\nu_\omega$,
\begin{equation}\label{eq:support-measure-invariance}
  g_*\Lambda_\omega^L=\Lambda_{\Theta\omega}^L.
\end{equation}
Consequently, for every such $(\omega,g)$,
\begin{equation}\label{eq:support-bundle-invariance}
 A(g)(\theta)L(\omega,\theta)
 =L(\Theta\omega,\theta+a(g))
 \qquad\text{for }m\text{-a.e. }\theta.
\end{equation}
The exceptional Haar-null set may depend on $(\omega,g)$.
\end{lemma}

\begin{proof}
The defining relation \eqref{eq:a.e.-bundle-invariance}, integrated against Haar measure, gives
\begin{equation}\label{eq:a.e.-measure-transport}
  g_*\Lambda_\omega^L=\Lambda_{\Theta\omega}^L
\end{equation}
for $P(\dd\omega)\nu_\omega(\dd g)$-almost every $(\omega,g)$.  To upgrade this relation, choose a countable convergence-determining family $(\varphi_\ell)_{\ell\geq1}$ in $C(Q_k)$.  For each $\ell$, define
\[
 D_\ell(\omega,g)
 :=\int_{Q_k}\varphi_\ell(g\cdot y)\dd\Lambda_\omega^L(y)
 -\int_{Q_k}\varphi_\ell(y)\dd\Lambda_{\Theta\omega}^L(y).
\]
Each $D_\ell$ is jointly measurable. Fubini's theorem and a
countable intersection give a full-measure Borel set $\Omega_L$
such that, for every $\omega\in\Omega_L$ and every $\ell$,
\[
  D_\ell(\omega,g)=0
  \qquad\text{for }\nu_\omega\text{-a.e. }g.
\]
For fixed $\omega$, the map $g\mapsto D_\ell(\omega,g)$ is continuous because the action \eqref{eq:grassmann-action} is continuous and $Q_k$ is compact.  Its zero set is closed and has full $\nu_\omega$-measure, hence contains $\supp\nu_\omega$.  Therefore all the functions $D_\ell(\omega,\cdot)$ vanish on $\supp\nu_\omega$.  Since the family $(\varphi_\ell)$ determines probability measures, \eqref{eq:support-measure-invariance} follows.

Both sides of \eqref{eq:support-measure-invariance} have Haar torus marginal.  Disintegrating over $\TT^d$, the left side has conditional measure at phase $\theta+a(g)$ equal to $\delta_{A(g)(\theta)L(\omega,\theta)}$, while the right side has conditional measure $\delta_{L(\Theta\omega,\theta+a(g))}$. Uniqueness of Rokhlin disintegration gives equality of these Dirac masses for Haar-almost every $\theta$, which is \eqref{eq:support-bundle-invariance}.
\end{proof}

\end{document}